\documentclass{amsart}
\usepackage{amsmath}
\usepackage{amsthm}
\usepackage{amssymb}
\usepackage{mathrsfs}
\usepackage[curve,arrow,matrix]{xy}
\theoremstyle{plain}
\newtheorem{thm}{Theorem}[section]
\newtheorem{prop}[thm]{Proposition}
\newtheorem{lemma}[thm]{Lemma}

\newtheorem{cor}[thm]{Corollary}

\newcommand{\deqno}{\refstepcounter{thm}(\thethm)}
\theoremstyle{definition}
\newtheorem{defn}[thm]{Definition}
\newtheorem{exam}[thm]{Example}

\theoremstyle{remark}
\newtheorem{remk}[thm]{Remark}

\newcounter{item}
\newenvironment{rlist}[1][1]{\begin{list}
  {\textup{(\roman{item})}}{\usecounter{item} \setcounter{item}{#1}\addtocounter{item}{-1}
  \setlength{\itemsep}{0ex}
  \setlength{\topsep}{0ex} \setlength{\parsep}{0ex} \setlength{\labelwidth}{15mm}
   \setlength{\leftmargin}{10mm} } }{\end{list}}
\newenvironment{alist}[1][1]{\begin{list}
  {\textup{(\alph{item})}}{\usecounter{item} \setcounter{item}{#1}\addtocounter{item}{-1}
  \setlength{\itemsep}{0ex}
  \setlength{\topsep}{0ex} \setlength{\parsep}{0ex} \setlength{\labelwidth}{15mm}
  \setlength{\leftmargin}{10mm} } }{\end{list}}
\newenvironment{nlist}[1][1]{\begin{list}
  {\textup{(\arabic{item})}}{\usecounter{item} \setcounter{item}{#1}\addtocounter{item}{-1}
  \setlength{\itemsep}{0ex}
  \setlength{\topsep}{0ex} \setlength{\parsep}{0ex} \setlength{\labelwidth}{15mm}
  \setlength{\leftmargin}{10mm} } }{\end{list}}

\newcommand{\bt}{{\bf t}}
\newcommand{\bx}{{\bf x}}
\newcommand{\by}{{\bf y}}

\newcommand{\F}{{\bf F}}
\newcommand{\Fe}{{\bf F}^e}
\newcommand{\nn}{\mathbb{N}}

\newcommand{\sB}{\mathscr{B}}

\newcommand{\fp}{\mathfrak{p}}
\newcommand{\fq}{\mathfrak{q}}

\newcommand{\tensor}{\otimes}
\newcommand{\Hom}{\textup{Hom}}

\newcommand{\Tor}{\textup{Tor}}

\newcommand{\Ann}{\textup{Ann}}

\newcommand{\hgt}{\textup{ht}}
\newcommand{\im}{\textup{Im}}

\newcommand{\Mod}{\textup{Mod}}

\newcommand{\ra}{\rightarrow}
\newcommand{\RA}{\Rightarrow}
\newcommand{\incl}{\hookrightarrow}

\newcommand{\isom}{\cong}

\newcommand{\wh}{\widehat}

\begin{document}

\title{Big Cohen-Macaulay Algebras and Seeds}
\author{Geoffrey D. Dietz}

\address{Department of Mathematics,
Gannon University, 109 University Square, Erie, Pennsylvania 16541}
\email{gdietz@member.ams.org}
\thanks{The author was supported in part by a VIGRE grant from
 the National Science Foundation.}
\thanks{Transactions of the American Mathematical Society, 359 (2007), No. 12,
5959--5989.}

%    General info

%\date{\today}


\begin{abstract}
In this article, we delve into
the properties possessed by algebras, which
we have termed \textit{seeds}, that map
to big Cohen-Macaulay algebras.
We will show that over a complete local domain
of positive characteristic any two big Cohen-Macaulay
algebras map to a common big Cohen-Macaulay
algebra. We will also strengthen Hochster and Huneke's ``weakly
functorial" existence result for big Cohen-Macaulay
algebras by showing that the
seed property is stable under base change between
complete local domains of positive characteristic.
We also show that every seed over a positive characteristic
ring $(R,m)$ maps to a balanced big Cohen-Macaulay $R$-algebra that
is an absolutely integrally closed, $m$-adically separated,
quasilocal domain.
\end{abstract}

\maketitle

%\tableofcontents

\section{Introduction}
\label{intro}

Methods used by M.\ Hochster and C.\ Huneke during their
development of tight closure led to the remarkable result that
$R^+$ (the integral closure of a domain $R$ in an algebraic
closure of its fraction field) is a \textit{balanced big Cohen-Macaulay
algebra} over $R$ when $R$ is an excellent local domain of
positive characteristic \cite{HH92}. A \textit{big Cohen-Macaulay algebra}
over a local ring $(R,m)$ is an $R$-algebra $B$ such that some system
of parameters of $R$ is a regular sequence on $B$. It is 
\textit{balanced} if \textit{every} system of parameters of $R$ is
a regular sequence on $B$. While Hochster had shown the
existence of big Cohen-Macaulay modules in equal characteristic
\cite{Ho75}, this new result was the first proof that big
Cohen-Macaulay algebras existed. Big Cohen-Macaulay algebras also
exist in equal characteristic 0 \cite{HH95} and
in mixed characteristic when $\dim R\leq 3$ \cite{Ho02}.
(The latter result follows from Heitmann's proof
of the direct summand conjecture for mixed characteristic
rings of dimension three \cite{Heit02}.)
Their existence is important as
it gives new proofs for many of the local homological conjectures,
such as the direct summand conjecture, monomial conjecture, and
vanishing conjecture for maps of $\Tor$. 

After defining \textit{seeds} as algebras over a
local ring $R$ that map to a (balanced) big Cohen-Macaulay
$R$-algebra and presenting some basic
properties, we characterize seeds in terms of the
existence of \textit{durable colon-killers}
(see Theorem~\ref{ckseed}) for positive characteristic
rings.
One of our most useful results is given
in Theorem~\ref{intseed}, where we
show that if $R$ is a local Noetherian
ring of positive characteristic, $S$ is
a seed over $R$, and $T$ is an
integral extension of $S$, then
$T$ is also a seed over $R$. We can
view this theorem as a generalization
of the existence of big Cohen-Macaulay
algebras over complete local domains
of positive characteristic; see Corollary~\ref{bigCMexist}.

We also define a class of \textit{minimal seeds},
which, analogously to the \textit{minimal solid
algebras} of Hochster (see \cite{Ho94}), are seeds that have no proper
homomorphic image that is also a seed. While
it has only been shown that Noetherian solid
algebras map onto minimal solid algebras in \cite{Ho94},
we have shown in Proposition~\ref{minseed} that
every seed maps onto a minimal seed. Furthermore,
like minimal solid algebras, we have shown
that minimal seeds are domains in positive characteristic;
see Proposition~\ref{mindom}. As a result, every seed
in positive characteristic over a local ring $(R,m)$
maps to a balanced big Cohen-Macaulay
algebra that is a quasilocal, absolutely integrally closed, $m$-adically
separated domain; see Theorem~\ref{nicebigCM}.

Among the most intriguing results of this article are 
Theorems~\ref{product} and \ref{basechange}. Both results concern seeds
over a complete local domain $R$ of positive characteristic.
Theorem~\ref{product} shows that the tensor product of two seeds
is still a seed, just as the tensor product of two solid modules
is still solid (see Proposition~\ref{solidprop} here). As an
immediate consequence, if $B$ and $B'$ are big Cohen-Macaulay
algebras over $R$, then $B$ and $B'$ both map to a common big
Cohen-Macaulay algebra $C$, which shows that the class of big
Cohen-Macaulay algebras over $R$ forms a directed system in this
sense.

Theorem~\ref{basechange} shows that if $R\ra S$ is a map of
positive characteristic complete local domains, and $T$ is a seed
over $R$, then $T\tensor_R S$ is a seed over $S$, so that the seed
property is stable under this manner of base change. This theorem
greatly strengthens the ``weakly functorial" existence result of
Hochster and Huneke in \cite[Theorem 3.9]{HH95}, where they show
that given complete local domains of equal characteristic $R\ra
S$, there exists a balanced big Cohen-Macaulay $R$-algebra $B$ and a 
balanced big
Cohen-Macaulay $S$-algebra $C$ such that $B\ra C$ extends the map
$R\ra S$. Our theorem shows that if we have $R\ra S$ in positive
characteristic, and $B$ is \textit{any} big Cohen-Macaulay
$R$-algebra, then there exists a big Cohen-Macaulay $S$-algebra
$C$ that fills in a commutative square:
$$
\xymatrix{ B\ar[r] & C\\ R\ar[r]\ar[u] & S\ar[u]}
$$

\section{Notation and Background}

All rings throughout are commutative with identity. All modules are unital.

\subsection{Big Cohen-Macaulay modules and algebras}
\begin{defn}
A \textit{big Cohen-Macaulay module} $M$ over a local Noetherian
ring $(R,m)$ is an $R$-module such that some system of
parameters for $R$ is a regular sequence on $M$. It is
\textit{balanced} if every system of parameters is
a regular sequence. (Since the
definition of a regular sequence $x_1,\ldots,x_d$ requires that
$(x_1,\ldots,x_d)M\neq M$, we have $mM\neq M$.) If $M=B$ is an
$R$-algebra, then $B$ is called a (balanced) \textit{big Cohen-Macaulay
algebra}.
\end{defn}

The terminology ``big" refers to the fact that $M$ is not
necessarily a finitely generated $R$-module. We say
that a partial system of parameters of $(R,m)$ is a
\textit{possibly improper regular sequence} on $M$ if all
relations on the parameters are trivial, but $mM\neq M$ does not
necessarily hold. Similarly, we can have \textit{possibly improper
big Cohen-Macaulay modules}.

The question ``When do big Cohen-Macaulay modules or algebras exist?" has important
applications in commutative algebra.
Over regular rings, the answer is simple.

\begin{prop}[p.77, \cite{HH92}]\label{regflatbigCM}
Let $R$ be a regular Noetherian ring, and let $M$ be an $R$-module.
Then $M$ is a balanced big Cohen-Macaulay module over $R$ if and only
if $M$ is faithfully flat over $R$.
\end{prop}

In \cite{Ho75}, Hochster showed that big Cohen-Macaulay modules
exist over all equicharacteristic local Noetherian rings. The
first significant existence proof of big Cohen-Macaulay algebras
came from a celebrated theorem of Hochster and Huneke. For a
domain $R$, let $R^+$ denote the integral closure of $R$ in an
algebraic closure of its fraction field. The ring $R^+$ is called
the \textit{absolute integral closure} of $R$ and is not
Noetherian in general.

\begin{thm}[Theorem 5.15, \cite{HH92}]\label{R+bigCM}
If $R$ is an excellent local domain, then $R^+$ is a balanced big
Cohen-Macaulay $R$-algebra.
\end{thm}

Using this result, Hochster and Huneke were also able to establish
a ``weakly functorial" existence of big Cohen-Macaulay algebras
over all equicharacteristic local rings. The term
\textit{permissible} used in the following theorem refers to a map
$R\ra S$ such that every minimal prime $Q$ of $\wh{S}$, with $\dim
\wh{S}/Q=\dim \wh{S}$, lies over a prime $P$ of $\wh{R}$ that
contains a minimal prime $\fp$ of $\wh{R}$ satisfying $\dim
\wh{R}/\fp=\dim \wh{R}$.

\begin{thm}[Theorem 3.9, \cite{HH95}]\label{wfbigCM}
We may assign to every equicharacteristic local ring $R$ a 
balanced big
Cohen-Macaulay $R$-algebra $\sB(R)$ in such a way that
if $R\ra S$ is a permissible local homomorphism of
equicharacteristic local rings, then we obtain a
homomorphism $\sB(R)\ra\sB(S)$ and a commutative diagram:
$$
\xymatrix{ \sB(R)\ar[r] & \sB(S) \\
R\ar[r]\ar[u] & S\ar[u] }
$$
\end{thm}

A key tool in the proof of this result is the construction
of big Cohen-Macaulay algebras using \textit{algebra
modifications}. Since we will make great use of algebra
modifications in this article, it will be helpful to review
some definitions and useful properties now.

\begin{defn}
Given a local Noetherian ring $R$, an $R$-algebra $S$,
and a relation $sx_{k+1} = \sum_{i=1}^k x_is_i$ in $S$,
where $x_1,\ldots,x_{k+1}$ is part of a system of parameters
of $R$, the $S$-algebra
$$
T:=\frac{S[U_1,\ldots,U_k]}{s-\sum_{i=1}^k x_iU_i}
\leqno\deqno\label{algmoddef}
$$
is called an \textit{algebra modification of $S$ over $R$}.
\end{defn}

We can also construct an $S$-algebra $\Mod(S/R) = \Mod_1(S/R)$
by adjoining infinitely many indeterminates and killing the appropriate
relations (as above) so that \textit{every} relation in $S$ on a
partial system of parameters from $R$ is trivialized in
$\Mod_1(S/R)$. Inductively define $\Mod_n(S/R)
 = \Mod(\Mod_{n-1}(S/R))$ and then define
$\Mod_\infty(S/R)$ as the direct limit of the
$\Mod_n(S/R)$. We have now formally trivialized all possible
relations on systems of parameters from $R$ and done so in
a way that is universal in the following sense.

\begin{prop}[Proposition 3.3b, \cite{HH95}]\label{modbigCM}
Let $S$ be an algebra over the local Noetherian ring $R$. Then
$\Mod_\infty(S/R)$ is a possibly improper balanced big Cohen-Macaulay
$R$-algebra. It is a proper balanced big Cohen-Macaulay algebra if and only if
$S$ maps to some balanced big Cohen-Macaulay $R$-algebra.
\end{prop}

While $\Mod_\infty(S/R)$ is a rather large and cumbersome object,
we can study it in terms of finite sequences of algebra modifications.
Given a Noetherian local ring $(R,m)$ and an $R$-algebra $S$, we
set $S^{(0)}:= S$ and then inductively define $S^{(i+1)}$ to be
an algebra modification of $S^{(i)}$ over $R$. We then obtain
a finite \textit{sequence of algebra modifications}
$$
S = S^{(0)}\ra S^{(1)}\ra \cdots \ra S^{(h)}
$$
for any $h\in \nn$. We call such a sequence \textit{bad} if
$mS^{(h)}=S^{(h)}$. We will frequently use the following
proposition.

\begin{prop}[Proposition 3.7, \cite{HH95}]\label{seqmod}
Let $R$ be a local Noetherian ring, and let $S$ be an $R$-algebra.
$\Mod_\infty(S/R)$ is a proper balanced big Cohen-Macaulay \mbox{$R$-algebra}
if and only if no finite sequence of algebra modifications
is bad.
\end{prop}

\subsection{The Frobenius endomorphism}

Throughout this article we will often work with rings of
positive prime characteristic $p$. We will always let $e$ denote
a nonnegative integer and let $q$ denote $p^e$, a power
of $p$. Thus, the phrase ``for all $q$" will mean ``for all powers
$q=p^e$ of $p$."

Every characteristic $p$ ring $R$ comes equipped with a
\textit{Frobenius endomorphism}
$F_R:R\ra R$, which maps
$r\mapsto r^p$. We can compose this map with itself to obtain the iterations
$F_R^e:R\ra R$,
which map $r\mapsto r^q$.
Associated to these maps
are the \textit{Peskine-Szpiro} (or \textit{Frobenius})
\textit{functors} $\Fe_R$.
If we let $S$ denote the ring
$R$ viewed as an $R$-module via the $e^{th}$-iterated Frobenius
endomorphism, then $\Fe_R$
is the covariant functor $S\tensor_R -$ which
takes $R$-modules to $S$-modules and so takes
$R$-modules to $R$-modules since $S=R$ as a ring. Specifically,
if $R^m\ra R^n$ is a map of free $R$-modules given by the matrix
$(r_{ij})$, then we may apply $\Fe_R$ to this map to obtain a map
between the same $R$-modules given by the matrix
$(r^q_{ij})$. For cyclic modules
$R/I$, $\Fe_R(R/I)=R/I^{[q]}$, where
$$
I^{[q]} := (a^q \,|\, a\in I)R
$$
is the \textit{$q^{th}$ Frobenius power} of the ideal $I$.
For modules $N\subseteq M$, we will denote the image of
$\Fe_R(N)$ in $\Fe_R(M)$ by $N_M^{[q]}$,
and we will denote
the image of $u\in N$ inside of $N_M^{[q]}$ by $u^q$.

\subsection{Tight closure}

The operation of tight closure
was developed by Hochster and Huneke in the late
1980s and early 1990s as a method for proving (often reproving with
dramatically shorter proofs) and generalizing theorems for
rings containing a positive characteristic field. See \cite{HH90}
for an introduction.

We will denote
the complement in $R$ of the set of minimal primes by
$R^\circ$.

\begin{defn} For a Noetherian ring $R$ of characteristic $p>0$
and finitely generated modules $N\subseteq M$, the \textit{tight
closure} $N_M^*$ of $N$ in $M$ is
$$
N_M^* := \{ u\in M\,|\, cu^q\in N_M^{[q]}\, \mbox{for all}\, q\gg 1, \,
\mbox{for some}\, c\in R^\circ\}.
$$
In the case that $M=R$ and $N=I$, $u\in I^*$
if and only if
there exists $c\in R^\circ$ such that $cu^q\in I^{[q]}$, for all
$q\gg 1$.
\end{defn}

A very powerful tool in the application of tight closure is
the notion of a \textit{test element}.

\begin{defn} For a Noetherian ring of positive characteristic
$p$, an element $c\in R^\circ$ is called a \textit{test element}
if for all ideals $I$ and all $u\in I^*$, we have
$cu^q\in I^{[q]}$ for all $q\geq 1$. If, furthermore, $c$ is a
test element in every localization of $R$, then
$c$ is \textit{locally stable}, and if it is also
a test element in the completion of
every localization, then $c$ is
\textit{completely stable}.
\end{defn}

Test elements exist in very general settings.

\begin{thm}[Theorem 6.1, \cite{HH94sm}]\label{testexist} If $R$ is
a reduced, excellent local ring of positive characteristic,
then $R$ has a completely stable test element.
\end{thm}

Below is a collection of important tight closure results that
we will refer to during this article. These properties help
make tight closure a ``good" closure operation for
solving problems related to the local homological
conjectures in positive characteristic.

\begin{thm}\label{TCprop}
Let $R$, $S$ be Noetherian rings of positive prime
characteristic $p$, and let $N\subseteq M$ be
finitely generated $R$-modules.
\begin{alist}
\item \cite[Theorem 4.4, Proposition 8.7]{HH90} If $R$ is regular, then $N_M^*=N$.

\item (\textbf{colon-capturing}) \cite[Theorem 4.7]{HH90}
Let $R$ be module-finite and torsion-free over a regular
domain $A$. Let $x_1,\ldots,x_n\in A$ be parameters in
$R$. Then $(x_1,\ldots,x_{n-1})R:_R x_n\subseteq ((x_1,\ldots,x_{n-1})R)^*$.

\item (\textbf{persistence}) \cite[Theorem 6.24]{HH94sm}
Let $R$, $S$ be excellent rings with $R\ra S$ a
homomorphism. If $w\in N_M^*$, then
$1\tensor w\in \im(S\tensor_R N\ra S\tensor_R M)_{S\tensor_R M}^*$.

\item \cite[Corollary 5.23]{HH94}
Let $S$ be a module-finite extension of $R$. If $1\tensor u\in
\im(S\tensor_{R} N\ra S\tensor_{R} M)_{S\tensor_{R} M}^*$,
calculated over $S$, then $u\in N_M^*$. In particular,
$(IS)^*_S\cap R \subseteq I_R^*$ for all ideals $I$ of $R$.
\end{alist}
\end{thm}

\subsection{Solid Algebras and Modules}

Hochster introduced the notion of \textit{solid} modules and algebras
in \cite{Ho94} in an attempt to define a characteristic free
notion of tight closure.

\begin{defn} If $R$ is a Noetherian domain, then an $R$-module
$M$ is \textit{solid} if $\Hom_R(M,R)\neq 0$. If $M=S$ is
an $R$-algebra, then $S$ is \textit{solid} over $R$
if it is solid as an $R$-module.
\end{defn}

The following properties have analogues in the world of
seeds and big Cohen-Macaulay modules. See Corollary~\ref{bigCMbasechange}
and Theorems~\ref{product} and \ref{basechange}.

\begin{prop}[Section 2, \cite{Ho94}]\label{solidprop}
Let $R$ be a Noetherian domain.
\begin{alist}
\item Let $S$ be a module-finite domain extension of $R$, and let
$M$ be an $S$-module. Then $M$ is solid over $S$ if and only if
$M$ is solid over $R$.

\item If $M$ and $N$ are solid $R$-modules, then $M\tensor_R N$ is solid.

\item (\textbf{persistence of solidity}) Let $R\ra S$ be any map
of Noetherian domains. If $M$ is a solid $R$-module,
then $S\tensor_R M$ is a solid $S$-module.
\end{alist}
\end{prop}

Whether seeds
and solid algebras are the same in positive characteristic
is an interesting open question. An answer to one direction
was given by Hochster.

\begin{thm}[Corollary 10.6, \cite{Ho94}]\label{bigCMsolid}
Let $R$ be a complete local domain. An $R$-algebra that
has an $R$-algebra map to a big Cohen-Macaulay algebra
over $R$ is solid.
\end{thm}

We also need a connection between tight closure and
big Cohen-Macaulay algebras.

\begin{thm}[Theorem 11.1, \cite{Ho94}]\label{bigCMTC}
Let $R$ be a complete local domain of positive characteristic, and
let $N\subseteq M$ be finitely generated $R$-modules. If $u\in M$,
then $u\in N_M^*$ if and only if there exists
a (balanced) big Cohen-Macaulay $R$-algebra $B$ such that
$1\tensor u\in \im(B\tensor_R N\ra B\tensor_R M)$.
\end{thm}

\section{Definition and Properties of Seeds}

\begin{defn} For any local Noetherian ring $(R,m)$, an $R$-algebra
$S$ is called a 
\textit{seed} over $R$ if $S$ maps to a big Cohen-Macaulay
$R$-algebra. \end{defn}

Bartijn and Strooker \cite[Theorem 1.7]{Bar-Str}
show that the $m$-adic separated completion 
of a big Cohen-Macaulay algebra (module) is 
a balanced big Cohen-Macaulay algebra (module).
The definition of seed is then unchanged by requiring
a balanced big Cohen-Macaulay algebra.
Using this terminology,
Proposition~\ref{seqmod} implies that $S$ is a seed
if and only if $S$ does not have a bad sequence of
algebra modifications. Based on this characterization of seeds,
a direct limit of seeds is still a seed.

\begin{lemma}\label{limit} Let $(R,m)$ be a local Noetherian ring, and
let $S=\varinjlim_\lambda S_\lambda$ be a direct limit of
a directed set of $R$-algebras. Then $S$ is a seed if and only if
each $S_\lambda$ is a seed. \end{lemma}
\begin{proof} Since $S$ is an $S_\lambda$-algebra for all $\lambda$,
if $S$ is a seed, then so is each $S_\lambda$. Conversely,
suppose that $S$ is not a seed. We will find an $S_\lambda$ that
also has a bad sequence of modifications.

As $S$ is not a seed, it has a bad sequence of algebra modifications.
It is then straightforward to show that this bad sequence
over $S$ induces
a bad sequence of modifications over some $S_\lambda$
because a finite sequence
involves only finitely many relevant elements. These elements
generate an $S_\lambda$ over $R$ that is not a seed.
\end{proof}

In positive characteristic, we can use
the Frobenius endomorphism and its iterates to map any seed to a
reduced and perfect balanced big Cohen-Macaulay algebra.
We will let $R^\infty$ denote the direct limit of
the directed system
$$R\ra \F(R)\ra\F^2(R)\ra\cdots\ra\Fe(R)\ra\cdots,$$
where $\Fe$ is the iterated Peskine-Szpiro functor. If $R$ is
reduced, then $R^\infty$ can be obtained by adjoining all $q^{th}$
roots to $R$. Notice that $R^\infty = (R_{\textup{red}})^\infty$,
where $R_{\textup{red}}$ is the quotient of $R$ obtained by
killing all nilpotents, and that $R^\infty$ is always reduced and
contains all of its $q^{th}$ roots.

\begin{lemma}\label{red/perfect} Let $R$ be a local Noetherian ring
of positive characteristic $p$. If $B$ is a balanced big Cohen-Macaulay $R$-algebra,
then $B^\infty$ is a reduced balanced big Cohen-Macaulay algebra containing all
of its $q^{th}$ roots. Moreover, if $B$ is quasilocal, then $B^\infty$ is also
quasilocal.
\end{lemma}
\begin{proof} A direct limit of balanced big Cohen-Macaulay algebras is a
balanced big Cohen-Macaulay algebra. If $B$ is quasilocal, then all maps in
the direct limit are local so that $B^\infty$ is quasilocal.
\end{proof}

We will show later (Proposition~\ref{domredmod})
that any reduced seed (in any
characteristic) can be modified into a reduced 
balanced big Cohen-Macaulay algebra.

\begin{lemma}\label{local} Let $(R,m)$ be a local Noetherian ring.
If $B$ is a (balanced) big Cohen-Macaulay $R$-algebra and $\fp$ is any prime ideal of $B$
containing $mB$, then $B_{\fp}$ is also a (balanced) big Cohen-Macaulay $R$-algebra.
Moreover, if $B$ is reduced (resp., $R$ has positive characteristic
and $B$ is reduced and perfect), then $B_{\fp}$ is still
reduced (resp., reduced and perfect). \end{lemma}
\begin{proof} Given $x_{k+1}(r/u)\in (x_{1},\ldots,x_{k})B_{\fp}$,
where $x_{1},\ldots,x_{k+1}$ is part of a system of parameters for
$R$, we may assume $u=1$ in showing that the relation is trivial.
Therefore, there
exists a $v\in B\setminus \fp$ such that $x_{k+1}(rv)\in
(x_{1},\ldots,x_{k})B$. Since $B$ is big Cohen-Macaulay, $rv\in
(x_{1},\ldots,x_{k})B$, and thus $r/1\in (x_{1},\ldots,x_{k})B_{\fp}$
as needed. Furthermore, since $mB\neq B$ and $\fp\supseteq mB$, we see
that $mB_{\fp}\neq B_{\fp}$, and so $B_{\fp}$ is a (balanced) 
big Cohen-Macaulay algebra.

The other claims follow from the following easy lemma.
\end{proof}

\begin{lemma}\label{localredperf} If $S$ is any reduced ring and $U$ is
a multiplicatively closed set in $S$, then $U^{-1}S$ is also reduced.
If, in addition, $S$ has positive characteristic and is perfect, then
$U^{-1}S$ is also perfect. \end{lemma}
\begin{proof} The first claim is well known.

Now suppose that $S$ has positive characteristic and is perfect.
Given $s/u$ in $U^{-1}S$, with $s\in S$ and $u\in U$, we have
$s=a^{q}$ and $u=b^{q}$, where $a$ and $b$ are in $S$. Then
$ab^{q-1}/u$ is an element of $U^{-1}S$ and is a $q^{th}$ root of
$s/u$.
\end{proof}

We can also use the separated
completion of a big Cohen-Macaulay $R$-algebra with respect to the maximal ideal
of $R$ to give us an $m$-adically
separated balanced big Cohen-Macaulay algebra while preserving the other properties we
have worked with earlier.

\begin{lemma}\label{complredperf}
If $A$ is a reduced and perfect ring of positive
characteristic $p$, and $I$ is an ideal of $A$, then the $I$-adic
completion $\wh{A}$ of $A$ is reduced and perfect.
\end{lemma}
\begin{proof} By definition,
$$\wh{A} = \{a=(a_1,a_2,a_3,\ldots)\in\prod_j A/I^j \, |\,
a_k\equiv a_j\pmod{I^j},\quad\forall k>j\}.$$
If $a^n =0$ in $\wh{A}$, then there exists $q$,
a power of $p$, such that $a^q = 0$,
so that $a_k^q\in I^k$ for all $k$. Given any index $j$, there exists an
integer $k(j)$ such
that $I^{k(j)}\subseteq (I^j)^{[q]}$. Therefore, for any $j$, we can find
$k(j)\geq j$ such that
$a_{k(j)}^q\in (I^j)^{[q]}$. Since $A$ is perfect, $a_{k(j)}\in I^j$, and
since $k(j)\geq j$,
we have $a_j\in I^j$. Hence, $a=0$, and $\wh{A}$ is reduced.

Given $a=(a_1,a_2,\ldots)\in\wh{A}$, we will now find an element
$b\in\wh{A}$ such that
$b^q=a$. Indeed, let $b=(a_{k(1)}^{1/q},a_{k(2)}^{1/q},\ldots)$,
where $k(j)$ is chosen
so that $k(j)\geq j$, $k(j)\geq k(j-1)$, and
$I^{k(j)}\subseteq (I^j)^{[q]}$. If $i\geq k(j)$, then $a_i\equiv a_{k(j)}
\pmod{I^{k(j)}}$, so that $a_i\equiv a_{k(j)} \pmod{(I^j)^{[q]}}$. Since
$A$ is perfect, we
can take $q^{th}$ roots to see that $a_i^{1/q}\equiv a_{k(j)}^{1/q}
\pmod{I^j}$, which
shows that $b$ is a well-defined element of $\wh{A}$. Finally, $b^q =
(a_{j(1)},a_{j(2)},\ldots)$, which is easily seen to be equal to $a$.
\end{proof}

This lemma is the last piece we need to show that seeds map to big
Cohen-Macaulay algebras with certain rather useful properties.
We
will show later that seeds map to big Cohen-Macaulay algebras
with even stronger properties. See Theorem~\ref{nicebigCM}.

\begin{prop}\label{goodbigC-M} Let $(R,m)$ be a Noetherian local
ring of positive characteristic. Every seed over $R$ maps to a 
balanced big
Cohen-Macaulay $R$-algebra $B$ that is reduced, perfect, quasilocal, and
$m$-adically separated. \end{prop}
\begin{proof} Use Lemmas~\ref{red/perfect}, \ref{local},
\ref{complredperf}, and \cite[Theorem 1.7]{Bar-Str}.
\end{proof}

\section{Colon-killers and Seeds}

Hochster and Huneke used colon-killers (also called
\textit{Cohen-Macaulay multipliers}) in \cite{HH92}
as tools for proving the existence of big Cohen-Macaulay algebras in
positive characteristic. Not surprisingly, the existence of
such elements in algebras over a local ring will be useful
in determining whether an algebra is a seed or not.
We shall work with a generalized version of their definition
and will define a special class of colon-killers that will help us
determine when an $R$-algebra is a seed.

\begin{defn}\label{ckdefn}
Let $R$ be a local Noetherian ring, $S$ an $R$-algebra,
and $M$ an arbitrary $S$-module.
An element $c\in S$ is a \textit{colon-killer} for $M$
over $R$ if
$$
c((x_1,\ldots,x_k)M:_M x_{k+1})\subseteq (x_1,\ldots,x_k)M,
$$
for each partial system of parameters $x_1,\ldots,x_{k+1}$
in $R$.
\end{defn}

We will soon prove that a colon-killer for an $S$-module
$M$ over $R$ has a power that is
a colon-killer for $M$ over $S$ when $S$ is an integral
extension of $R$. First, we need the next lemma
connecting colon-killers and Koszul homology.

\begin{lemma}\label{ckKoszul} Let $R$ be a local Noetherian ring,
let $S$ be an $R$-algebra, and
let $M$ be an arbitrary $S$-module. If $c\in S$ is nonzero, then
the following are equivalent:
\begin{rlist}
\item Some power of $c$ is a colon-killer for $M$ over $R$.
\item Some power of $c$ kills the Koszul homology
modules $H_i(x_1,\ldots,x_k;M)$ for all $i\geq 1$ and
all partial systems of parameters $x_1,\ldots,x_k$.
\item Some power of $c$ kills $H_1(x_1,\ldots,x_k;M)$
for all partial systems of parameters $x_1,\ldots,x_k$.
\end{rlist}
\end{lemma}
\begin{proof} (ii) $\RA$ (iii) is obvious. For (iii) $\RA$ (i), let
$\bx = x_1,\ldots, x_{k}$ be part of a system of
parameters for $R$, and let $\bx' = x_1,\ldots,x_{k-1}$. We obtain
a short exact sequence
$$
0\ra \frac{H_i(\bx';M)}{x_k H_i(\bx';M)}\ra H_i(\bx;M)\ra
\Ann_{H_{i-1}(\bx';M)}x_k \ra 0 \leqno\deqno\label{koszulses}
$$
for all $i$ from \cite[Corollary 1.6.13(a)]{BH}. In the case $i=1$,
we see that there is a surjection of $H_1(\bx;M)$ onto the module
$((\bx')M:_M x_k)/(\bx')M$, which implies that the latter
module is killed by the same power of $c$ that
kills the former.

For (i) $\RA$ (ii), assume without loss of generality that
$c$ itself is a colon-killer for $M$ over $R$. We will use induction
on $k$ to show that $c^{2^{k-1}}$ kills
$H_i(x_1,\ldots,x_k;M)$ for $i\geq 1$.
If $k=1$, then $H_1(x_1;M)$ is the only
nonzero Koszul homology module, and it is isomorphic
to $\Ann_M x_1 = (0:_M x_1)$. Since $c$ is a
colon-killer for $M$, $c$ kills $H_1(x_1;M)$.
Now let $k\geq 2$, $\bx=x_1,\ldots,x_k$, $\bx'=x_1,\ldots,x_{k-1}$,
and suppose that $c^{2^{k-2}}$ kills $H_i(\bx';M)$ for
$i\geq 1$. Using the sequence (\ref{koszulses}),
we see that $c^{2^{k-1}}$ kills $H_i(\bx;M)$ for
all $i\geq 2$ by the inductive hypothesis, and
$c^{2^{k-2}+1}$ kills $H_1(\bx;M)$ by the
inductive hypothesis together with
$c$ being a colon-killer. Therefore,
if $N=2^{\dim R -1}$, then $c^N$ kills
all of the relevant Koszul homology modules.
\end{proof}

From this lemma, we can obtain our result
on colon-killers.

\begin{prop}\label{ckbasechange}
Let $S$ be a Noetherian local ring that
is an integral extension of a
local Noetherian ring $R$. Let $M$ be an arbitrary $S$-module.
If $c\in R$ kills the Koszul homology modules
$H_i(x_1,\ldots,x_k;M)$ for all $i\geq 1$ and
all partial systems of parameters $x_1,\ldots,x_k$ in $R$,
then $c^N$ kills $H_i(y_1,\ldots,y_k;M)$ for all $i\geq 1$ and
all partial systems of parameters $y_1,\ldots,y_k$ in $S$,
for some $N$. Consequently, if $c$ is a colon-killer
for $M$ over $R$, then
a power of $c$ is a colon-killer for $M$ over $S$.
\end{prop}
\begin{proof} Based on the previous lemma, it is enough
to prove the first claim. Let $\by = y_1,\ldots, y_k$ be
part of a system of parameters for $S$.
Since $R\incl S$ is integral, $R/((\by)\cap R)\incl S/(\by)$
is also integral, so that
$$
\dim R/((\by)\cap R) = \dim S/(\by)
= \dim S - k = \dim R - k.
$$
We claim that $(\by)\cap R$ contains a partial system
of parameters $\bx = x_1,\ldots,x_k$ for $R$. Indeed, we proceed
by induction on $k$, where the case $k=0$ is trivial. We
can then assume without loss of generality that $k=1$ and
obtain our result from the general fact that if $I$ is an ideal of $R$
such that $\dim R/I < \dim R$, then $I$ contains a parameter of $R$.
If not, then for every $x\in I$, there exists a prime ideal $\fp$ of
$R$ such that $x\in\fp$ and $\dim R/\fp =\dim R$. Therefore,
$I$ is contained in the union of such prime ideals, and so by prime
avoidance, $I$ is contained in a prime $\fp$ such that $\dim R/\fp =
\dim R$. We have a contradiction, which implies that $I$ contains
a parameter.

Thus, there exists part of a system of parameters
$\bx=x_1,\ldots,x_k$ of $R$
such that $(\bx)S\subseteq (\by)$. Then a power of $c$ kills
$H_i(\by;M)$ for all $i\geq 1$ by the following lemma, where
the power depends only on $k$, which in turn is bounded by
$\dim R$.
\end{proof}

\begin{lemma} Let $S$ be any ring and $M$ any $S$-module.
Suppose $(x_1,\ldots,x_k)S\subseteq (y_1,\ldots,y_k)S=(\by)S$
and that
$c\in S$ kills $H_i(x_{1},\ldots,x_{m};M)$ for all $i\geq 1$ and
all $1\leq m\leq k$. Then
$c^{D(m)}$ kills $H_{k+1-m}(\by;M)$ for $1\leq m\leq k$,
where $D(1)=1$, and $D(m)=2^{k-2}D(m-1)+2$ for $m\geq 2$.
\end{lemma}
\begin{proof}
We will use induction on $m$.
For $m=1$, we need $c$ to kill $H_{k}(\by;M)$, but
$$
H_{k}(\by;M)\isom \Ann_{M}(\by)\subseteq \Ann_{M}(x_{1})\isom
H_{1}(x_{1};M),
$$
which implies what we want.

Now, let $m\geq 2$.
The hypothesis on $c$ together with Lemma~\ref{ckKoszul} implies
that $c$ is a colon-killer for $M$ with respect to subsequences
of $x_{1},\ldots,x_{k}$, but this implies that $c$ is a colon-killer
for $M/x_{1}M$ with respect to subsequences of $x_{2},\ldots,x_{k}$.
The proof of Lemma~\ref{ckKoszul} then implies that $c^{2^{k-2}}$
kills $H_{i}(x_{2},\ldots,x_{m};M/x_{1}M)$ for all $i\geq 1$ and
all $2\leq m\leq k$.

For the induction, suppose that $c^{D(m-1)}$ kills $H_{k+2-m}(\by;M)$.
Consider the exact sequence
$$
0\ra \Ann_{M}x_{1}\ra M\overset{x_{1}}\ra M\ra M/x_{1}M\ra 0,
$$
from which we obtain two short exact sequences
$$
\begin{array}{c}
0\ra \Ann_{M}x_{1}\ra M \ra x_{1}M\ra 0, \\
0\ra x_{1}M\ra M\ra M/x_{1}M\ra 0.
\end{array}
$$
These sequences then induce long exact sequences in Koszul homology:
$$
\begin{array}{c}
H_{i+1}(\by;M)\overset{f_{i+1}}\ra H_{i+1}(\by;x_{1}M)\ra
H_{i}(\by;\Ann_{M}x_{1})\ra H_{i}(\by;M)\overset{f_{i}}\ra
H_{i}(\by;x_{1}M), \\

H_{i+1}(\by;x_{1}M)\overset{g_{i+1}}\ra H_{i+1}(\by;M)\ra
H_{i+1}(\by;\frac{M}{x_{1}M})\ra H_{i}(\by;x_{1}M)\overset{g_{i}}\ra H_{i}(\by;M),
\end{array}
$$
which in turn yield short exact sequences for all $i\geq 0$:
$$
0\ra \frac{H_{i+1}(\by;x_{1}M)}{f_{i+1}(H_{i+1}(\by;M))}\ra
H_{i}(\by;\Ann_{M}x_{1})\ra \ker(f_{i})\ra 0, \leqno(*_{i})
$$
$$
0\ra \frac{H_{i+1}(\by;M)}{g_{i+1}(H_{i+1}(\by;x_{1}M))}\ra
H_{i+1}(\by;M/x_{1}M)\ra \ker(g_{i})\ra 0. \leqno(\#_{i})
$$

Since the map of homology
$$
H_{i}(\by;M)\overset{f_{i}}\ra H_{i}(\by;x_{1}M)
\overset{g_{i}}\ra H_{i}(\by;M)
$$
is induced by the composition $M\ra x_{1}M\ra M$, which
is multiplication by $x_{1}$, and since $x_{1}\in(\by)S$,
$g_{i}\circ f_{i} = 0$ for all $i\geq 0$.

Since $c\Ann_{M}x_{1} = 0$ by hypothesis, $(*_{i})$ implies
that
$$
cH_{i+1}(\by;x_{1}M)\subseteq f_{i+1}(H_{i+1}(\by;M))\subseteq
\ker(g_{i+1}),
$$
for all $i\geq 0$. Applying the inductive hypothesis to $M/x_{1}M$
implies that $c^{2^{k-2}D(m-1)}$ kills $H_{k+2-m}(\by;M/x_{1}M)$.
Then $(\#_{k+1-m})$ shows that $c^{2^{k-2}D(m-1)}
\ker(g_{k+1-m}) = 0$. Therefore, $c^{2^{k-2}D(m-1)+1}$ kills
$H_{k+1-m}(\by;x_{1}M)$.
To finish, notice
that
$$
c^{2^{k-2}D(m-1)+1}=c^{D(m)-1},
$$
which kills
the image of $H_{k+1-m}(\by;M)$ inside $H_{k+1-m}(\by;x_{1}M)$
under the map $f_{k+1-m}$. Thus, $c^{D(m)-1}H_{k+1-m}(\by;M)$
is contained in $\ker(f_{k+1-m})$, which is killed by $c$, using
$(*_{k+1-m})$, and so $c^{D(m)}$ kills $H_{k+1-m}(\by;M)$, as
needed.
\end{proof}

Since 1 is a colon-killer for a balanced big Cohen-Macaulay module,
the previous result gives us a result about
big Cohen-Macaulay modules and base change.

\begin{cor}\label{bigCMbasechange}
Let $S$ be a Noetherian local ring that
is also an integral extension of a
local Noetherian ring $R$. If $M$ is an $S$-module and
a balanced big Cohen-Macaulay $R$-module,
then $M$ is a balanced big Cohen-Macaulay $S$-module.
\end{cor}

We will now
introduce another notion of a colon-killer that
will be very useful for us in the following sections when we need
to determine whether a ring is a seed.

\begin{defn} For a local Noetherian ring $(R,m)$ and an $R$-algebra $S$,
an element $c\in S$ is called a \textit{weak durable colon-killer} over $R$
if for some system of parameters $x_1,\ldots,x_n$ of $R$,
$$c( (x_1^{t},\ldots,x_k^{t})S:_S x_{k+1}^{t} ) \subseteq
(x_1^{t},\ldots,x_k^{t})S,$$
for all $1\leq k\leq n-1$ and all $t\in\nn$, and if for any $N\geq 1$,
there exists $k\geq 1$ such that $c^N\not\in m^k S$. An element $c\in
S$ will be simply called a \textit{durable colon-killer} over $R$
if it is a weak durable colon-killer for every system of parameters of $R$.\end{defn}

Notice that if $S=R$, then all colon-killers in $R$ that are not
nilpotent are durable colon-killers. So, if $R$ is reduced, all
nonzero colon-killers are durable colon-killers. We can now use
the existence of durable colon-killers to characterize when an
algebra is a seed by adapting the proof of \cite[Theorem
11.1]{Ho94}.

\begin{thm}\label{ckseed} Let $(R,m)$ be a local Noetherian ring of
positive characteristic $p$, and let
$S$ be an $R$-algebra. Then $S$ is a seed if and only if there is a map
$S\ra T$ such that $T$ has a (weak) durable colon-killer $c$.
\end{thm}
\begin{proof} If $S$ is a seed, then $S\ra B$, for some balanced
big Cohen-Macaulay $R$-algebra
$B$. As pointed out above, 1 is a durable colon-killer in $B$, so $T=B$
will suffice.
For the converse, we will modify the proof of
\cite[Theorem 11.1]{Ho94} to obtain our result. We will show that the
existence of a (weak) durable colon-killer in an $S$-algebra $T$ implies
that $S$ is a seed. (All parenthetical remarks will apply to the case
that $T$ only possesses a weak durable colon-killer.)

Suppose that $S\ra T$ is a map
such that $T$ has a (weak) durable colon-killer $c$
(with respect to a fixed system of parameters in $R$).
Let $S^{(0)}:=S$, and given $S^{(i)}$ for $0\leq i\leq t-1$, let
$$S^{(i+1)} := \frac{S^{(i)}[U_1^{(i)},\ldots, U_{k_i}^{(i)}]}
{s^{(i)} - \sum_{j=1}^{k_i} x_j^{(i)}U_j^{(i)}},$$
where $x_1^{(i)},\ldots,x_{k_i+1}^{(i)}$
is a system of parameters for $R$ 
(the fixed system of parameters in the latter case), and
$x_{k_i+1}^{(i)}s^{(i)} = \sum_{j=1}^{k_i} x_j^{(i)}s_j^{(i)}$
is a relation in $S^{(i)}$. Then
$$S=S^{(0)}\ra S^{(1)}\ra S^{(2)}\ra \cdots \ra S^{(t)}$$
is a finite sequence of algebra modifications. Suppose to the contrary
that $S$ is not a seed and the sequence is bad, so that $1\in mS^{(t)}$.
(In the weak
case, we are supposing that $S$ does not map to an $S$-algebra where
the fixed system of parameters is a regular sequence.) We can then
write
$$1=\sum_{j=1}^{n}r_{j}w_{j}, \leqno\deqno\label{badeqn}$$
where $r_{j}\in m$ and $w_{j}\in S^{(t)}$ for all $j$.

We will construct inductively homomorphisms
$\psi_e^{(i)}$ from each $S^{(i)}$ to $\Fe(T)$ forming a commutative diagram:
$$\xymatrix{\Fe(T)_{c}\ar@{=}[r] & \Fe(T)_{c}\ar@{=}[r] & \Fe(T)_{c}\ar@{=}[r] &
\cdots\ar@{=}[r] & \Fe(T)_{c} \\
S^{(0)}\ar[r]\ar[u]^{\psi_e^{(0)}} & S^{(1)}\ar[r]\ar[u]^{\psi_e^{(1)}} &
S^{(2)}\ar[r]\ar[u]^{\psi_e^{(2)}} & \cdots\ar[r] & S^{(t)}.\ar[u]^{\psi_e^{(t)}}  }
$$

In order to construct the maps we need to keep track of
bounds, independent of $q=p^e$, associated with the images of
certain elements of each $S^{(i)}$. For all $1\leq i\leq t$, we will use
reverse induction to define a finite subset $\Gamma_i$ of $S^{(i)}$
and positive integers $b(i)$. We will then inductively define positive
integers $\beta(i)$ and $B(i)$, which will be the necessary bounds.

First, let $$\Gamma_{t}:=\{w_{1},\ldots,w_{n}\},$$
where the $w_{j}$ are from relation (\ref{badeqn}).
Now, given $\Gamma_{i+1}$ (with $0\leq i\leq t-1$), each element can be written
as a polynomial in the $U_j^{(i)}$ with coefficients
in $S^{(i)}$. Let $b(i+1)$ be the largest degree of any such polynomial.
For $i\geq 1$, let
$\Gamma_{i}$ be the set of all coefficients of these polynomials
together with
$s^{(i)}, s_1^{(i)},\ldots,s_{k_{i}}^{(i)}$.
Now define $\beta(1):=1$, $B(1):=b(1)$, and given $B(i)$ for $1\leq i\leq t-1$,
let
$$\beta(i+1):= B(i)+1 \qquad \textup{and} \qquad
B(i+1) := \beta(i+1)b(i+1)+B(i).$$
Notice that, as claimed, all
$\beta(i)$ and $B(i)$ are independent of $q$.

Fix $q=p^{e}$. By hypothesis, we have a map $S^{(0)}=S\ra T$ that can
be naturally
extended to a map $\psi_{e}^{(0)}:S^{(0)}\ra \Fe(T)_{c}$ by composing
with $T\ra \Fe(T)_c$. We next
define a map $\psi_{e}^{(1)}:S^{(1)}\ra \Fe(T)_{c}$ that extends
$\psi_{e}^{(0)}$, maps the $U_{j}^{(0)}$ to the cyclic
$\Fe(T)$-submodule $c^{-1}\Fe(T)=c^{-\beta(1)}\Fe(T)$ in $\Fe(T)_{c}$, and
maps $\Gamma_{1}$ to $c^{-b(1)}\Fe(T)=c^{-B(1)}\Fe(T)$ inside
$\Fe(T)_{c}$. To do this we need only find appropriate images of
the $U_{j}^{(0)}$ such that the image of $s^{(0)} -
\sum_{j=1}^{k_0} x_j^{(0)}U_j^{(0)}$ maps to 0.

Since
$$x_{k_0+1}^{(0)}s^{(0)} = \sum_{j=1}^{k_0} x_j^{(0)}s_j^{(0)},$$
we have
$$(x_{k_0+1}^{(0)})^{q}\psi_{e}^{(0)}(s^{(0)}) = \sum_{j=1}^{k_0}
(x_j^{(0)})^{q}\psi_{e}^{(0)}(s_j^{(0)})$$
in $\Fe(T)_{c}$, where
the image of $\psi_{e}^{(0)}$ is contained in the image of
$\Fe(T)$ inside $\Fe(T)_{c}$. As $c$ is a (weak) durable colon-killer in
$T$, and so also a (weak) colon-killer in $\Fe(T)$,
$$c\psi_{e}^{(0)}(s^{(0)}) = \sum_{j=1}^{k_0}
(x_j^{(0)})^{q}\sigma_{j}^{(0)},$$
where the $\sigma_{j}^{(0)}$ are in the image of $\Fe(T)$ in
$\Fe(T)_{c}$. If we define $\psi_{e}^{(1)}$ such that
$U_{j}^{(0)}\mapsto c^{-1}\sigma_{j}^{(0)}$, then we have
accomplished our goal for $\psi_{e}^{(1)}$ because the elements of
$\Gamma_{1}$ can be written as polynomials in the $U_{j}^{(0)}$
of degree at most $b(1)=B(1)$ with
coefficients in $S$.

Now suppose that for some $1\leq i\leq t-1$ we have a map
$\psi_e^{(i)}: S^{(i)}\ra \Fe(T)_{c}$,
where the $U_j^{(i-1)}$ all map to $c^{-\beta(i)}\Fe(T)$, and
$\Gamma_i$ maps to $c^{-B(i)}\Fe(T)$. We will extend $\psi_e^{(i)}$
to a map from $S^{(i+1)}$ such that
each $U_j^{(i)}$ maps to $c^{-\beta(i+1)}\Fe(T)$, and $\Gamma_{i+1}$ maps
to $c^{-B(i+1)}\Fe(T)$.

In order to simplify notation, we drop many of the $(i)$ labels on
parameters. Then
$$S^{(i+1)}=\frac{S^{(i)}[U_{1},\ldots,U_{k}]}{s-\sum_{j=1}^{k}x_{j}U_{j}}.$$
Since $s$ and the $s_{j}$ (in the
relation $x_{k+1}s = \sum_{j=1}^{k}x_{j}s_{j}$ in $S^{(i)}$) are in
$\Gamma_{i}$, we can write
$$\psi_{e}^{(i)}(s) = c^{-B(i)}\sigma \qquad
\mbox{and}\qquad \psi_{e}^{(i)}(s_{j}) = c^{-B(i)}\sigma_{j},
\leqno\deqno $$
where $\sigma$
and the $\sigma_{j}$ are elements in the image of $\Fe(T)$ in
$\Fe(T)_{c}$. Hence,
$$x_{k+1}^{q}\psi_{e}^{(i)}(s) =
\sum_{j=1}^{k}x_{j}^{q}\psi_{e}^{(i)}(s_{j})$$
in $\Fe(T)_{c}$. Multiplying through by $c^{B(i)}$ yields
$$x_{k+1}^{q}\sigma = \sum_{j=1}^{k}x_{j}^{q}\sigma_{j}$$
in the image of $\Fe(T)$ in $\Fe(T)_{c}$.
Using our (weak) colon-killer $c$, we have
$$c\sigma = \sum_{j=1}^{k}x_{j}^{q}\tau_{j},$$
where $\tau_j$ is an element in the image of $\Fe(T)$ in
$\Fe(T)_{c}$. Therefore,
$$\psi_{e}^{(i)}(s) = \sum_{j=1}^{k}x_{j}^{q}(c^{-B(i)-1}\tau_{j})$$
in $\Fe(T)_{c}$.

We now have a well-defined
map $\psi_{e}^{(i+1)}:S^{(i+1)}\ra \Fe(T)_{c}$ extending
$\psi_{e}^{(i)}$ given by
$$
\psi_{e}^{(i+1)}(U_{j}) = c^{-B(i)-1}\tau_{t} = c^{-\beta(i+1)}\tau_t
$$
such that the $U_{j}$ map to
$c^{-\beta(i+1)}\Fe(T)$, and $\Gamma_{i+1}$ maps to
$c^{-B(i+1)}\Fe(T)$ since
$$
B(i+1) = \beta(i+1)b(i+1)+B(i)
$$
and these elements can be written as
polynomials in the $U_{j}$ of degree at most $b(i+1)$ with
coefficients in $\Gamma_{i}$.

We can finally conclude that, for all $q=p^{e}$, there exists a map
$$\psi_{e}^{(t)}:S^{(t)}\ra \Fe(T)_{c}$$
such that the equation (\ref{badeqn}) that puts $1\in mS^{(t)}$ maps to
$$1 = \sum_{j=1}^{n}r_{j}^{q}\psi_{e}^{(t)}(w_{j}).$$
If we let $B:=B(t)$, then each $\psi_{e}^{(t)}(w_{j})$ is in
$c^{-B}\Fe(T)$ as $\Gamma_{t}$ contains these elements. Multiplying
through by $c^{B}$, we see that
$c^{B} \in m^{[q]}T$ for all $q\geq 1$, where $B$ is independent
of $q$, which implies that $c^B\in m^kT$ for all $k\geq 1$. Since
$c$ is a (weak) durable colon-killer, we have a contradiction. Therefore,
no such finite sequence of modifications of $S$ is bad.

In the case of the durable colon-killer, we see that $S$ 
maps to a balanced big Cohen-Macaulay algebra over $R$.
In the weak case, $S$ maps to a big Cohen-Macaulay algebra that
can be made balanced by
\cite[Theorem 1.7]{Bar-Str}. In either event, $S$ is a seed.
\end{proof}

Later we will use this result as a piece of the proof that
integral extensions of seeds are seeds. We will also use durable
colon-killers to obtain our results concerning
when the seed property is preserved by base change.

\section{Minimal Seeds}

\begin{defn} For a Noetherian local ring $(R,m)$, an $R$-algebra
$S$ is a \textit{minimal seed} if $S$ is a seed over $R$ and no proper
homomorphic image of $S$ is a seed over $R$. \end{defn}

\begin{exam} \begin{nlist}
\item If $R$ is a Cohen-Macaulay ring, then $R$ is a minimal seed.

\item If $R$ is an excellent local domain of positive prime characteristic, 
then $R^+$ is a balanced big Cohen-Macaulay algebra over $R$ and a minimal seed.
\end{nlist}
\end{exam}

We also point out the following easy, but useful, fact about minimal
seeds.

\begin{lemma}\label{mininj}
A seed $S$ over a local Noetherian ring $R$ is a minimal seed if and only if every map
from $S$ to a (balanced) big Cohen-Macaulay $R$-algebra $B$ is injective if and only if
every map to a seed over $R$ is injective.
\end{lemma}

A very important question about minimal seeds is whether
or not every seed maps to a minimal seed.

\begin{prop}\label{minseed} Let $R$ be a local Noetherian ring, and let $S$ be a
seed over $R$. Then $S/I$ is a minimal seed for some ideal $I\subseteq S$.
\end{prop}
\begin{proof} Let $\Sigma$ be the set of all ideals $J$ of $S$ such that
$S/J$ is a seed. If $\Sigma$ contains a maximal element $I$, then
$S/I$ will be a minimal seed. Let $J_1\subset J_2\subset \cdots$
be a chain of ideals in $\Sigma$, and let $J=\bigcup_{k} J_k$.
Then $S/J = \varinjlim_k S/J_k$, and since each $S/J_k$ is a seed, 
Lemma~\ref{limit} implies that $S/J$ is a seed as well. By Zorn's Lemma,
$\Sigma$ has a maximal element $I$.
\end{proof}

Now that we know minimal seeds exist, we would like to know whether they are
domains or not.
After dealing with integral extensions of seeds, we will
prove in Section~\ref{moreseedprops} that in positive characteristic
minimal seeds are domains. In the meantime, we will point out
that minimal seeds are reduced in positive characteristic.

\begin{prop} \label{minred}
Let $S$ be a minimal seed over a local ring $R$ of positive characteristic. Then
$S$ is a reduced ring. \end{prop}
\begin{proof} By Lemma~\ref{red/perfect}, there exists a reduced big Cohen-Macaulay
algebra $B$ such that $S\ra B$. Since $S$ is minimal, Lemma~\ref{mininj}
implies that $S\incl B$. As $B$ has no nilpotents, neither does $S$.
\end{proof}

\section{Integral Extensions of Seeds}

The primary goal of this section is to prove that integral
extensions of seeds are seeds in positive characteristic.
Since all integral extensions are direct limits of module-finite
extensions, with Lemma~\ref{limit} we can concentrate on module-finite
extensions of seeds.
We begin the argument by proving that the problem can be
reduced to
a much more specific problem, which we attack
by constructing a durable colon-killer in a certain module-finite extension
of a big Cohen-Macaulay algebra. Our first reduction of the problem will be that
we can assume we are considering a module-finite extension of
a balanced big Cohen-Macaulay algebra that is reduced, quasilocal, and
$m$-adically separated.

\begin{lemma}\label{reduction} Let $(R,m)$ be a local Noetherian ring
of positive characteristic,
and let $S$ be a seed with $T$ a module-finite extension of $S$. Suppose
that any module-finite extension of a reduced, quasilocal, $m$-adically
separated balanced big Cohen-Macaulay algebra is a seed. Then $T$ is a seed.
\end{lemma}
\begin{proof} By Propositions~\ref{minseed} and \ref{minred},
$S/I$ is a minimal reduced seed for some ideal $I$. Since $S/I$ is
reduced, $I$ is a radical ideal and so is an integrally closed ideal.
Therefore, $IT\cap S = I$ so that $S/I$ injects into $T/IT$, which is thus
a module-finite extension of $S/I$. Since $T$ is a seed if
$T/IT$ is a seed, we can now assume that $S$ is a reduced seed.

By Lemma~\ref{mininj} and Proposition~\ref{goodbigC-M},
there exists a commutative square
$$\xymatrix{B\ar[r] & C \\
S\,\ar@{^{(}->}[u]\ar@{^{(}->}[r] & T,\ar[u]}$$ where $B$ is a
reduced, quasilocal, and $m$-adically separated balanced big Cohen-Macaulay
algebra, and $C:=T\tensor_S B$. If we can show that the upper map
$B\ra C$ is injective, then $C$ will be a module-finite extension
of $B$, and our hypotheses will imply that $C$ and $T$ are seeds.
Since $B$ is reduced, the next lemma allows us to reach our goal.
\end{proof}

\begin{lemma} If $S$ is a ring, $T$ is an integral extension of $S$, and
$B$ is a reduced extension of $S$, then $B$ injects into
$C:=T\tensor_S B$.
\end{lemma}
\begin{proof} We will first prove the claim in the case that $S$, $T$, and
$B$ are all domains. Let $K$ be the algebraic closure of the
fraction field of $S$, and let $L$ be the algebraic closure of the
fraction field of $B$. Since $T$ is an integral extension domain
of $S$, we have the following diagram:
$$
\xymatrix{ B\,\ar@{^{(}->}[rr] && L \\
S\,\ar@{^{(}->}[u]\ar@{^{(}->}[r] & T\,\ar@{^{(}->}[r] &
K.\ar@{^{(}->}[u]}
$$
Under the injection $K\incl L$, the ring $T$ maps
isomorphically to some subring $T'$ of $L$. Now let
$C'$ be the $S$-subalgebra of $L$ generated by $B$ and $T'$.
Since $C=T\tensor_S B$, we have a (surjective) map
$C\ra C'$ and a diagram
$$\xymatrix{ && C' \\
B\ar[r]\ar@/^/ @{^{(}->}[rru] & C\ar[ru] & \\
S\,\ar@{^{(}->}[u]\ar@{^{(}->}[r] & T\ar[u]\ar@/_/ @{^{(}->}[uur]
& }$$ Since $B$ injects into $C'$, $B$ also injects into $C$.

For the general case, since $B$ is reduced, it will suffice to
show that the kernel of $B\ra C$ is contained in every prime ideal
of $B$. Let $\fp$ be a prime of $B$, and let $\fp_0 := \fp\cap S$.
Since $T$ is integral over $S$, there exists a prime $\fq_0$ of
$T$ lying over $\fp_0$. If we put $D := T/\fq_0\tensor_{S/\fp_0}
B/\fp$, then we obtain the following commutative diagram:
$$\xymatrix{ & T/\fq_0\ar[rr] && D \\
T\ar[ur]\ar[rr] & \ar[u] & C\ar[ur] & \\
& S/\fp_0\ar@{^{(}-}[r]\ar@{^{(}-}[u] & \ar[r] & B/\fp\ar[uu] \\
S\,\ar[ur]\ar@{^{(}->}[rr]\ar@{^{(}->}[uu] && B\ar[uu]\ar[ur] &  }
$$

Since $T/\fq_0$ is still an integral extension of $S/\fp_0$ and $B/\fp$ is
a domain extension of $S/\fp_0$, the domain case of the proof shows
that $B/\fp$ injects into $D$. Therefore, if $b$ is in the kernel of
$B\ra C$, then $b$ is in the kernel of $B\ra D$, which implies
that $b\in\fp$, as desired.
\end{proof}

To finish our argument that module-finite
extensions of seeds are seeds, we will show that
we can extend the map $B\incl C$ to another
module-finite extension $B^{\#}\incl C^{\#}$ such that
$B^{\#}$ is a reduced, quasilocal, $m$-adically
separated balanced big Cohen-Macaulay algebra. Our new rings
will also have a nonzero
element $b\in B^{\#}$ such that $b$ multiplies $C^{\#}$ into
a finitely generated free $B^{\#}$-submodule of $C^{\#}$.
We will then finally show that $b$ is a durable colon-killer
in $C^{\#}$ so that Theorem~\ref{ckseed} implies that
$C$ is a seed.

We start the process by showing we can adjoin indeterminates
and then localize our ring $B$ without losing any of its
key properties.

\begin{defn} Let $(B,\fp)$ be a quasilocal ring. If $n,s\in\nn$
and $\{t_{ij}\, |\, i\leq n, j\leq s\}$ is a set of indeterminates
over $B$, then
$$
B^{\#(n,s)}:= B[t_{ij}\, |\, i\leq n,j\leq s]_{\fp B[t_{ij}]}.
$$
\end{defn}

\begin{lemma}\label{genmod} Let $(B,\fp)$ be a reduced, quasilocal,
$m$-adically separated balanced big Cohen-Macaulay $R$-algebra, where
$(R,m)$ is a local Noetherian ring. If $n,s\in\nn$, then
$B^{\#(n,s)}$ is also a reduced, quasilocal, $m$-adically
separated balanced big Cohen-Macaulay $R$-algebra.
\end{lemma}
\begin{proof} For the duration of the proof, $n$ and $s$ will be fixed, so
we will simply write $B^{\#}$
instead of $B^{\#(n,s)}$.
We will let $\bt$ denote the set of all $t_{ij}$.

As $B$ is reduced, $B[\bt ]$ is reduced
after adjoining indeterminates. By Lemma~\ref{localredperf},
$B^{\#}$ will also be reduced. As $\fp$ is prime in $B$,
the extension of $\fp$ to $B[\bt ]$ is also prime so that
it makes sense to localize at this ideal and end up with
$B^{\#}$ quasilocal.

The construction of $B^{\#}$ implies that $B^{\#}$ is a
faithfully flat extension of $B$. Therefore, for any ideals
$I$ and $J$ of $B$, we have $IB^{\#}:_{B^{\#}}JB^{\#} =
(I:_B J)B^{\#}$ (see \cite[Theorem 18.1]{N72}). This fact implies
that every system of parameters of $R$ will be a possibly improper
regular sequence on $B^{\#}$, because $B$ is a balanced big
Cohen-Macaulay algebra. Moreover, the faithful flatness also implies
that $mB^{\#}\neq B^{\#}$ as $mB\neq B$. Hence, $B^{\#}$
is a balanced big Cohen-Macaulay $R$-algebra.

To show that $B^{\#}$ is also $m$-adically separated will
take a little bit more effort. Suppose that an element
$F\in B^{\#}$ is in $m^NB^{\#}$ for every $N$. Multiplying through
by its denominator, we obtain such an element from $B[\bt ]$, so that
we may assume without loss of generality that $F$ is a polynomial
in $B[\bt ]$. Thus, for every $N$ there exists $G_N\not\in \fp B[\bt ]$
such that $G_N F\in m^NB[\bt ]$. It suffices to show that any polynomial
$G\not\in \fp B[\bt ]$ is not a zerodivisor modulo $m^N B[\bt ]$ for
any $N$. If we put $D:=B/m^N B$, then the image $\overline{G}$
of $G$ in $D[\bt ]$ is a polynomial not in $\fp D[\bt ]$, i.e.,
$\overline{G}$ is a polynomial whose coefficients generate
the unit ideal of $D$. To finish, it suffices to apply the following
general lemma.
\end{proof}

\begin{lemma} If $D$ is any ring and $G$ is a polynomial in
$D[t_\lambda \, |\, \lambda\in\Lambda]$ such that the coefficients of $G$
generate the unit ideal in $D$, then $G$ is not a zerodivisor.
\end{lemma}
\begin{proof} The lemma reduces to the Noetherian case, where it
follows from Corollary~2 on p.\ 152 of \cite{Mat}.
\end{proof}

\begin{lemma}\label{sharpcompat} If $(B,\fp)$ is a quasilocal
ring, and $B^{\#(n,s)}=B[t_{ij}\, |\, i\leq n,j\leq s]_{\fp B[t_{ij}]}$
for some $n$ and some $s$, then for any $k\leq n$
$$
B^{\#(n,s)}\isom (B^{\#(k,s)})^{\#(n-k,s)}.
$$
\end{lemma}
\begin{proof}
Let $\bx$ denote the indeterminates $t_{ij}$ such that
$1\leq i\leq k$ and $1\leq j\leq s$, and let $\by$ denote
the remaining indeterminates $t_{ij}$. Let $C :=
B^{\#(k,s)} = B[\bx]_{\fp B[\bx]}$, $Q$ be the maximal
ideal of $C$, and $D:= (B^{\#(k,s)})^{\#(n-k,s)} =C[\by]_{Q C[\by]}$.

Since $D$ is a $B$-algebra, there exists a unique
ring homomorphism $B[\bx,\by]\ra D$ that maps
the indeterminates $t_{ij}$ to their natural images
in $D$. We claim that the units in $B[\bx,\by]$
map to units in $D$ under this map. Indeed,
if $f(\bx,\by)$ is a unit in $B[\bx,\by]$, then
$f$ has some coefficient that is in $B\setminus \fp$.
If we rewrite
$$
f(\bx,\by) = g_k(\bx)\by^k +\cdots + g_1(\bx)\by +g_0(\bx),
$$
then some $g_i(\bx)$ is in $C\setminus Q$ so that
the image of $f$ is not in the maximal ideal of $D$.
We therefore have a ring homomorphism
$\phi:B^{\#(n,s)}\ra D$. Using the
previous lemma, it is easy to verify that
$\phi$ is injective. It is also routine
to check that $\phi$ is surjective.
\end{proof}

Using the construction
$B^{\#(n,s)}=B[t_{ij}\, |\, i\leq n,j\leq s]_{\fp B[t_{ij}]}$,
where $(B,\fp)$ is a quasilocal ring,
we also define
$$
^{\#(n,s)}M:=B^{\#(n,s)}\tensor_{B} M
$$
for any $B$-module $M$.
Since $B^{\#(n,s)}$ is faithfully flat over $B$, when $M=C$ is a
module-finite extension of $B$, we also have that $^{\#(n,s)}C$ is a
module-finite extension of $B^{\#(n,s)}$.

\begin{lemma}\label{modlemma} Let $B$ be a reduced quasilocal ring,
and let $M$
be a $B$-module generated by $m_{1},\ldots,m_{s}$
in $M$. Then there exists $k\leq s$ such that
$b(\,^{\#(k,s)}M)\subseteq G$, where $b$ is a nonzero
element of $B^{\#(k,s)}$ and $G$ is a finitely generated free
$B^{\#(k,s)}$-submodule of $^{\#(k,s)}M$.
\end{lemma}
\begin{proof} Throughout the proof, define
$$g_i:=t_{i1}m_1+\cdots + t_{is}m_s$$
in any $B^{\#(n,s)}$, where $i\leq n$.
Note that there exists a maximum $0\leq n\leq s$
such that the set $\{g_1,\ldots, g_n\}$
generates a $B^{\#(n,s)}$-free submodule of $^{\#(n,s)}M$, where
$B^{\#(0,s)}=B$, since $^{\#(n,s)}M$ has $s$ generators.
If $\alpha = (t_{ij})_{1\leq i,j\leq s}$, then $\det(\alpha)$ is not
in the unique maximal ideal of $B^{\#(s,s)}$, so that
$\alpha$ is an invertible matrix. As
$m_1,\ldots,m_s$ generate $^{\#(s,s)}M$ over $B^{\#(s,s)}$ and $\alpha$
is invertible, $g_1,\ldots,g_s$ also generate
$^{\#(s,s)}M$. If the $g_i$ are linearly independent over $B^{\#(s,s)}$, then
$^{\#(s,s)}M$ is a free module, and we can use $k=s$, $b=1$, and
$G= \,^{\#(s,s)}M$ to fulfill our claim.

Otherwise, the maximum value $n$ is strictly less than $s$,
and we put $k:=n+1$.
In this case, there exists a nonzero $b'\in B^{\#(k,s)}$ such that
$$
b'g_k \in (g_1,\ldots,g_{k-1})(^{\#(k,s)}M).
$$
Indeed, since
$n$ was chosen to be a maximum and $k=n+1$, there must
be a nontrivial relation $b'g_k = b_1g_1+\cdots +b_{k-1}g_{k-1}$ in
$^{\#(k,s)}M$. If $b'=0$, then we have a nontrivial relation
on $g_1,\ldots,g_{k-1}$. As $B^{\#(k,s)}\isom (B^{\#(k-1,s)})^{\#(1,s)}$
(by Lemma~\ref{sharpcompat}),
we see that $B^{\#(k,s)}$ is faithfully flat over $B^{\#(k-1,s)}$, and so
the nontrivial relation on $g_1,\ldots,g_{k-1}$ in $^{\#(k,s)}M$ implies
that there is a nontrivial relation on $g_1,\ldots,g_{k-1}$ in
$^{\#(k-1,s)}M$, a contradiction. Hence, $b'$ is nonzero as claimed.
Notice that the same argument implies that $g_1,\ldots,g_{k-1}$
still generate a finitely generated free submodule $G$ of
$^{\#(k,s)}M$.

We now claim that there exists a nonzero
$b\in B^{\#(k,s)}$ such that $b(^{\#(k,s)}M)\subseteq G$. If we
put $M_0:=(^{\#(k,s)}M)/G$,
and replace $m_1,\ldots,m_s$ and $g_k$ by their images in
$M_0$, then $b'$ kills $g_k$.
We intend to show that the annihilator of $M_0$ cannot be zero.
Suppose to the contrary that no nonzero element
of $B^{\#(k,s)}$ kills $M_0$. Then we have
an injective map
$B^{\#(k,s)}\incl M_0^{\oplus s}$ defined by
$b\mapsto (bm_1,\ldots,bm_s)$.

After clearing the denominator
on $b'$, we may assume that $b'$ is a polynomial
in $B[\bt]$, where $\bt$ denotes the set of all $t_{ij}$ with
$i\leq k$ and $j\leq s$. Write $b' = \sum_{\nu} c_{\nu} \bt^\nu$, where
$\nu$ is an $n\times s$ matrix of integers, and each $c_\nu$ is
in $B$.  Let $A_0$ be the prime ring of $B$, and let
$A$ be the $A_0$-subalgebra of $B$ (finitely) generated
by the nonzero $c_\nu$. Then $A$ is Noetherian and reduced.
If we let $\fq$ be the contraction of $\fp$ to $A$ and replace
$A$ by the local ring $A_\fq$, then $(A,\fq)$ is a reduced local Noetherian
subring of $B$. We then obtain an injective map
$A^{\#}\incl B^{\#(k,s)}$, where $A^{\#}$ denotes $A^{\#(k,s)}$.

Since $b'$ is in $A^{\#}$ (and still nonzero), we can define an $A^{\#}$-module
by
$$N:= \frac{A^{\#}m_1\oplus\cdots\oplus A^{\#}m_s}
{b'(t_{k1}m_1+\cdots+t_{ks}m_s)}.$$
There is then a natural map $N\ra M_0$ that induces a commutative
square:
$$\xymatrix{B^{\#(k,s)}\ar@{^{(}->}[r] & M_0^{\oplus s} \\
A^{\#}\ar@{^{(}->}[u]\ar[r] & N^{\oplus s}.\ar[u]}
$$
This implies that the map $A^{\#}\ra N^{\oplus s}$ defined
by $a\mapsto (am_1,\ldots,am_s)$ is injective (which also shows that
$N\neq 0$). Therefore, we may now
assume without loss of generality that $B$ is a reduced, local Noetherian
ring.

As $B$ is reduced and Noetherian, $B$ has finitely many minimal primes
$Q_1,\ldots, Q_h$ such that $\bigcap_i Q_i = 0$. Since
$b'$ is a nonzero polynomial in $B^{\#(k,s)}$, some coefficient of $b'$ is
not in some minimal prime $Q$. Thus, the image of $b'$ is still
nonzero in $(B_Q)^{\#(k,s)}$. Moreover, if $M_0':=(B_Q)^{\#(k,s)}\tensor_{B^{\#(k,s)}}
M_0$, then $M_0'$ is a finitely generated $(B_Q)^{\#(k,s)}$-module with
$b'(t_{k1}m_1+\cdots+t_{ks}m_s)=0$ and with an injection
$(B_Q)^{\#(k,s)}\incl (M_0')^{\oplus s}$ since
$(B_Q)^{\#(k,s)}\isom (B^{\#(k,s)})_{QB^{\#(k,s)}}$. (Again, this fact implies
that $M_0'$ is nonzero.) Since
$B$ is reduced and $Q$ is minimal, $B_Q$ is a field, and so
we can now assume that $B=K$ is a field.

In this final case, $K^{\#(k,s)}\isom K(\bt)$, and $M_0$ is a nonzero
module over a field so that $M_0$ is a nonzero free $K^{\#(k,s)}$-module.
Therefore, if $b'(t_{k1}m_1+\cdots+t_{ks}m_s)=0$ in $M_0$, then
$t_{k1}m_1+\cdots+t_{ks}m_s=0$ in $M_0$. This is impossible, however,
since the $t_{ij}$ are algebraically independent.

The resulting contradiction implies that $M_0$ is killed by some
nonzero element $b\in B^{\#(k,s)}$ in our original set-up, and
since $M_0$ was originally $(^{\#(k,s)}M)/G$, where $G$ is free
over $B^{\#(k,s)}$, the proof is complete.
\end{proof}

We are now ready to show that module-finite
extensions of sufficiently nice big Cohen-Macaulay algebras are indeed seeds. As
mentioned above, the key fact will be that the element $b$
constructed above is a durable colon-killer.

\begin{lemma}\label{mfck} Let $(R,m)$ be a local Noetherian ring of positive
characteristic, and let
$B$ be a reduced, quasilocal, $m$-adically separated balanced big Cohen-Macaulay
$R$-algebra. If $C$ is a module-finite extension of $B$, then $C$ is
a seed.
\end{lemma}
\begin{proof} By Lemma~\ref{genmod} and
the remarks made before the previous lemma, for any $k$,
we have a commutative square:
$$\xymatrix{ B^{\#(k,s)}\ar@{^{(}->}[r] & ^{\#(k,s)}C \\
B\,\ar[u]\ar@{^{(}->}[r] & C,\ar[u]}
$$
where the top map is also a module-finite extension of a reduced, quasilocal,
\mbox{$m$-adically} separated balanced big Cohen-Macaulay $R$-algebra, 
and $C$ is generated by $s$
elements as a $B$-module. After applying the previous lemma,
we may assume that there exists a nonzero element
$b\in B$ such that $b$ multiplies $C$
into a finitely generated free $B$-submodule $G$. In order to see that $C$ is
a seed, we show that $b$ is a durable colon-killer in $C$.

Indeed, let $x_1,\ldots, x_{t+1}$ be part of a system of parameters in $R$ and suppose
that $u\in (x_1,\ldots,x_t)C:_C x_{t+1}$. Then by construction
$bux_{t+1}\in (x_1,\ldots,x_t)G$, so as an element of $G$ we have
$(bu)\in (x_1,\ldots,x_t)G:_G x_{t+1}$. Since $B$ is a 
balanced big Cohen-Macaulay $R$-algebra and since
$G$ is a free $B$-module, $G$ is a balanced big Cohen-Macaulay $R$-module. Hence,
$bu\in (x_1,\ldots,x_t)G\subseteq (x_1,\ldots,x_t)C$ as $G$ is a submodule
of $C$.

Now, if $b^N\in \bigcap_t m^t C$ for some $N$, then
$b^{N+1}\in \bigcap_t m^t G$. Since $B$ is $m$-adically separated,
$\bigcap_t m^t B = 0$, and since $G$ is free over $B$, we also have
$\bigcap_t m^t G = 0$. As $G$ is a submodule of $C$ and the map
$B\ra C$ is an injection, $b^{N+1}=0$ in $B$, but $B$ reduced implies
that $b=0$, a contradiction. Therefore, $b$ is a durable colon-killer,
and $C$ is a seed by Theorem~\ref{ckseed}.
\end{proof}

We have now gathered together all of the tools that we will need to
prove the primary result of this section.

\begin{thm}\label{intseed} Let $(R,m)$ be a local Noetherian ring of
positive characteristic. If $S$ is a seed
and $T$ is an integral extension of $S$, then $T$ is a seed.
\end{thm}
\begin{proof} By Lemma~\ref{limit}, we may assume that $T$
is a module-finite extension of $S$ because integral
extensions are direct limits of module-finite extensions.
By Lemma~\ref{reduction}, we may assume that
$S=B$ is a reduced, quasilocal, $m$-adically
separated balanced big Cohen-Macaulay algebra. 
Finally, Lemma~\ref{mfck} implies
that $T$ is a seed.
\end{proof}

\begin{remk}
We feel it is worthwhile to point out that the hypothesis that our
base ring has positive characteristic is only required in two places:
(1) the existence of a durable colon-killer implies that an algebra is a
seed, and (2) all seeds map to a reduced big Cohen-Macaulay algebra. These facts
have proofs that rely heavily on the use of the Frobenius endomorphism, but
are the only two that we
can prove only in positive characteristic.
\end{remk}

We can also view the above theorem
as a generalization of the existence of big Cohen-Macaulay
algebras over complete local domains of positive
characteristic.

\begin{cor}[Hochster-Huneke]\label{bigCMexist} If $R$ is a complete
local domain of positive characteristic, then
there exists a balanced big Cohen-Macaulay algebra $B$
over $R$. \end{cor}
\begin{proof}
By the Cohen structure theorem, $R$ is a module-finite
extension of a regular local ring $A$. Since $A$ is
clearly a seed over itself, Theorem~\ref{intseed} implies
that $R$ is a seed over $A$ as well. Let $B$
be a balanced big Cohen-Macaulay algebra over $A$ such that
$B$ is also an $R$-algebra. By Corollary~\ref{bigCMbasechange}, 
$B$ is also a balanced big Cohen-Macaulay algebra over $R$.
\end{proof}

\section{More Properties of Seeds}
\label{moreseedprops}

In this section, we will show that all seeds in positive characteristic
can be mapped to quasilocal balanced big Cohen-Macaulay algebra domains
that are absolutely integrally closed and $m$-adically separated.
We start off the section by
delivering the promised proof that minimal seeds are domains. First,
we show that we can reduce to the case of a finitely generated minimal
seed.

\begin{lemma} Let $R$ be a local Noetherian ring. If all
finite type minimal seeds are domains, then
all minimal seeds are domains. \end{lemma}
\begin{proof} Let $S$ be an arbitrary minimal seed over $R$. Then
$S=\varinjlim_{\lambda\in\Lambda} S_\lambda$, where $\Lambda$
indexes the set of all finitely generated subalgebras of $S$.
Suppose that $S$ is not a domain and that $ab=0$ in $S$ with
$a,b\neq 0$. Since $S$ is a minimal seed, $S/aS$ and $S/bS$ are
not seeds. Let $\Lambda(a)$ and $\Lambda(b)$ be the subsets of
$\Lambda$ indexing all finitely generated subalgebras of $S$ that
contain $a$ and $b$, respectively. Then $S/aS =
\varinjlim_{\lambda\in \Lambda(a)} S_\lambda/aS_\lambda$, with a
similar result for $S/bS$. Since $S/aS$ and $S/bS$ are not seeds,
Lemma~\ref{limit} implies that there exists an $S_\alpha$
containing $a$ and an $S_\beta$ containing $b$ such that
$S_\alpha/aS_\alpha$ and $S_\beta/bS_\beta$ are not seeds.
Therefore, there exists a common $S_\gamma$ containing $a$ and $b$
such that $S_\gamma$ is not a seed modulo $aS_\gamma$ nor modulo
$bS_\gamma$. We also have $ab=0$ in $S_\gamma$ since $ab=0$ in
$S$. As $S$ is a seed, $S_\gamma$ is a seed. Since $S_\gamma$ is
also finitely generated as an $R$-algebra, $S_\gamma$  maps onto a
finitely generated minimal seed $T$. Therefore, $ab=0$ in $T$, and
as $T$ is a domain by hypothesis, $a=0$ or $b=0$ in $T$. Suppose
without loss of generality that $a=0$. This implies that the map
$S_\gamma \ra T$ factors through $S_\gamma/aS_\gamma$, which is
not a seed and so cannot map to any seed. We have a contradiction,
and so $S$ is a domain after all.
\end{proof}

\begin{prop}\label{mindom}
Let $R$ be a local Noetherian ring of positive characteristic
$p$. If $S$ is a minimal seed over $R$, then $S$ is a domain.
\end{prop}
\begin{proof} By the previous lemma, we can assume that $S$ is
finitely generated over $R$. By
Proposition~\ref{minred}, $S$ is Noetherian and reduced.
Let $\overline{S}$ be
the normalization of $S$ in its total quotient ring. Then
$\overline{S}$ is a finite direct product
of normal domains by Serre's Criterion
(see \cite[Theorem 11.5]{E}) and is an integral extension of $S$.
By our main result of the last
section, Theorem~\ref{intseed}, $\overline{S}$ is also a seed.
Since $\overline{S}$ is a seed and a finite product $D_{1}\times\cdots\times
D_{t}$ of domains, we will be done
once we have proven the following lemma.
\end{proof}

\begin{lemma} Let $(R,m)$ be a Noetherian local ring.
If $S=S_{1}\times\cdots\times S_{t}$, then $S$ is a seed over $R$
if and only if $S_{i}$ is a seed over $R$, for some $i$. \end{lemma}
\begin{proof} Clearly, if some $S_{i}$ is a seed, then $S$ is also a
seed. Suppose then that $S$ is a seed, but no $S_{i}$ is a seed. Since
$S$ is a direct product, if $S\ra B$, a balanced big Cohen-Macaulay algebra, then
$B\isom B_{1}\times\cdots\times B_{t}$, where each $B_{i}$ is an
$S_{i}$-algebra. We first claim that each $B_{i}$ is a possibly
improper balanced big Cohen-Macaulay algebra. 
Indeed, let $x_{k+1}b = \sum_{j=1}^{k}
x_{j}b_{j}$ be a relation in $B_{i}$ on a partial system of parameters
$x_{1},\ldots,x_{k+1}$ from $R$, and let $e_{i}$ be the idempotent
associated to $B_{i}$ in $B$. Therefore, $x_{k+1}(be_{i}) = \sum_{j=1}^{k}
x_{j}(b_{j}e_{i})$ is a relation in $B$, and so $be_{i} =
\sum_{j=1}^{k}x_{j}c_{j}$, for elements $c_{j}$ in $B$. Multiplying
this equation by $e_{i}$ yields $be_{i} =
\sum_{j=1}^{k}x_{j}(c_{j}e_{i})$ since $e_{i}^{2}=e_{i}$. If we let
$c_{j}'$ be the image of $c_{j}e_{i}$ in $B_{i}$, for all $j$, then
$b = \sum_{j=1}^{k}x_{j}c_{j}'$ in $B_{i}$, as claimed.

Now, if, as assumed, each $S_{i}$ is not a seed, then $1\in mB_{i}$,
for all $i$. Thus $e_{i}\in mB$, for all $i$, and so
$1=\sum_{i}e_{i}\in mB$, a contradiction. Therefore, some $S_{i}$ must
be a seed if $S$ is a seed.
\end{proof}

As a corollary, we will prove that each seed maps to a big Cohen-Macaulay algebra
that is also a domain. We must first
show that a domain seed can be modified into a balanced big Cohen-Macaulay algebra
domain. We will use the ``$\mathfrak{A}$-transform"
$$\Theta = \Theta(x,y;S) := \{u\in S_{xy} \,|\, (xy)^Nu\subseteq \im(S\ra S_{xy}),
\textup{ for some } N\}.$$
See \cite[Chapter V]{N65} and \cite[Section 12]{Ho94} for an
introduction to the properties of $\Theta$. The most useful property for
us is that if $x,y$ form part of a system of parameters in a local
Noetherian ring $R$ and $S$ is a seed over $R$, then $x,y$ become
a regular sequence on $\Theta$ (see \cite[Lemma 12.4]{Ho94}).
As a result any map from $S$ to
a balanced big Cohen-Macaulay $R$-algebra $B$ factors through $\Theta$.

Suppose now that $S$ is a seed over a local Noetherian ring $R$, and let
$$T = \frac{S[U_1,\ldots,U_k]}{(s-s_1U_1-\cdots -s_kU_k)}$$
be an algebra modification of $S$,
where $x_{k+1}s=x_1s_1+\cdots +x_ks_k$ is a relation in $S$ on
a partial system of parameters $x_1,\ldots,x_{k+1}$ in $R$. When
$k\geq 2$, we also have an induced relation on $x_1,\cdots,x_{k+1}$ in $\Theta =
\Theta(x_1,x_2;S)$ so that
$$T' = \frac{\Theta[U_1,\ldots,U_k]}{(s-s_1U_1-\cdots -s_kU_k)}$$
is an algebra modification of $\Theta$ over $R$. We will call $T'$ an
\textit{enhanced algebra modification of $S$ over $R$
induced by the relation $x_{k+1}s=x_1s_1+\cdots +x_ks_k$}.

Since any map from $S$ to a balanced big Cohen-Macaulay algebra $B$ factors through
$\Theta$, we obtain a commutative diagram
$$\xymatrix{ \Theta\ar[r] & T'\ar[dr] & \\
S\ar[u]\ar[r] & T\ar[u]\ar[r] & B}
$$
which shows that maps from algebra modifications of seeds to balanced big
Cohen-Macaulay algebras factor through the enhanced modification of $S$
when $T$ is a modification with respect to a relation on 3 or more parameters from $R$.
With this factorization, we can adapt the process
described in \cite[Section 3]{HH95}
to construct a balanced big Cohen-Macaulay algebra from a given seed as a very large
direct limit of enhanced and ordinary modifications.

\begin{lemma}\label{domredenmod} Let $(R,m)$ be a local Noetherian ring, and
let $S$ be a domain (resp., reduced). If $T$ is an enhanced algebra
modification of $S$ over $R$, then $T$ is also a domain (resp., reduced).
\end{lemma}
\begin{proof} Let $T$ be induced by the relation
$x_{k+1}s=x_1s_1+\cdots +x_ks_k$ in $S$, where $k\geq 2$. Then
$x_1,x_2$ forms a possibly improper regular sequence on
$\Theta=\Theta(x_1,x_2;S)$ (see \cite[Lemma 12.4]{Ho94}), and
so $x_1,s-x_1U_1-\cdots -x_kU_k$ is a possibly improper
regular sequence on
$\Theta[U_1,\ldots,U_k]$. (Any polynomial $f(U)$ that kills
$s-x_1U_1-\cdots -x_kU_k$ modulo $x_1$ has a highest degree term as a
polynomial in $U_2$, but this term is killed by $x_2$ modulo $x_1$.
Since $x_2$ is not a zerodivisor modulo $x_1$, the term must
be divisible by $x_1$. Hence, $f(U)$ is divisible by $x_1$,
and $s-x_1U_1-\cdots -x_kU_k$
is not a zerodivisor modulo $x_1$.)
Therefore, $x_1$ is not a zerodivisor on
$$
\Theta[U_1,\ldots,U_k]/(s-x_1U_1-\cdots -x_kU_k).
$$
The result now follows from the following short lemma.
\end{proof}

\begin{lemma} Let $A$ be a domain (resp., reduced). If $a$
and $x$ are elements of  $A$ such that $x$ is not a zerodivisor on
$A':=A[U]/(a-xU)$, then $A'$ is also a domain (resp.,
reduced).
\end{lemma}
\begin{proof} Since $x$ is not a zerodivisor on $A'$, we have an
inclusion $A'\incl (A')_x$ so that it suffices to show $(A')_x$
is a domain (resp., reduced). It is, however, easy to verify that
$(A')_x\isom A_x$ via the map that sends $U$ to $a/x$ (even without
any hypotheses on $A$ or $x$). Since
$A$ is a domain (resp., reduced), so is $A_x$.
\end{proof}

As promised, we now prove that one can modify a domain or reduced
seed into a balanced big Cohen-Macaulay algebra with the same property. 
Notice that the result is characteristic free.

\begin{prop}\label{domredmod} Let $R$ be a local Noetherian ring.
If $S$ is a seed and a domain (resp., reduced), then $S$ maps to a 
balanced big
Cohen-Macaulay algebra that is a domain (resp., reduced).
\end{prop}
\begin{proof} Any element of $S$ killed by a
parameter of $R$ will be in the kernel of any map to a balanced 
big Cohen-Macaulay
$R$-algebra $B$, so that the map $S\ra B$ factors through the
quotient of $S$ modulo the ideal of elements killed by a parameter
of $R$. When $S$ is a domain, this ideal is the zero ideal, and when
$S$ is reduced, this ideal is radical. Hence the quotient
is still a domain (resp., reduced). Without loss of generality,
we may then assume that no element of $S$ is killed by a parameter of
$R$. Thus, any nontrivial relation
$x_{k+1}s=x_1s_1+\cdots +x_ks_k$ in $S$ on part of a system
of parameters $x_1,\ldots,x_{k+1}$ from $R$ will have
$k\geq 1$.

If $k=1$, then an algebra modification with respect
to that relation factors through the $\mathfrak{A}$-transform
$\Theta(x_1,x_2;S)$. If $k\geq 2$,
then we can factor any algebra modification through
an enhanced algebra modification. Therefore,
we can map $S$ to a balanced big Cohen-Macaulay $R$-algebra that is constructed
as a very large direct limit of sequences of enhanced
algebra modifications and $\mathfrak{A}$-transforms.

Since $S$ is a domain (resp., a reduced ring) and since
the
\mbox{$\mathfrak{A}$-transform} of a domain (resp., reduced ring) is a
domain (resp., reduced),
Lemma~\ref{domredenmod} implies that all enhanced algebra
modifications and $\mathfrak{A}$-transforms in any sequence
will continue to be domains (resp., reduced).
Hence, $S$ maps to a balanced big Cohen-Macaulay algebra
$B$ that is a direct limit of domains (resp., reduced rings),
and so $B$ itself is a domain (resp., reduced).
\end{proof}

\begin{cor}\label{dom} Let $R$ be a local Noetherian ring
of positive characteristic. If $S$ is a seed, then $S$ maps to a 
balanced big
Cohen-Macaulay algebra domain. \end{cor}
\begin{proof} By Proposition~\ref{minseed}, $S$ maps to a minimal seed,
and Proposition~\ref{mindom} implies that the minimal seed is
a domain. The previous lemma then implies that a minimal seed can be
mapped to a balanced big Cohen-Macaulay algebra that is a domain.
\end{proof}

We can also show that all seeds map
to big Cohen-Macaulay algebras with a host of nice properties. We will use
an uncountable limit ordinal number in the proof and refer the reader to
\cite[Chapters 7 and 8]{HJ} for the properties
of such ordinal numbers.

\begin{thm}\label{nicebigCM}
Let $(R,m)$ be a local Noetherian ring of positive
characteristic. If $S$ is a
seed, then $S$ maps to an absolutely integrally closed, $m$-adically separated,
quasilocal balanced big Cohen-Macaulay algebra domain $B$.
\end{thm}
\begin{proof} We will construct $B$ as a direct limit of seeds indexed
by an uncountable ordinal number. Let $\beta$ be an uncountable
initial ordinal of cardinality $\aleph_1$.
Using transfinite induction, we will define an $S$-algebra
$S_\alpha$, for each ordinal number $\alpha<\beta$, and then
we will define $B$ to be the direct limit of all such $S_\alpha$.

Let $S_0=S$. Given a seed $S_\alpha$, we can form a sequence
$$
S_\alpha \ra S_\alpha^{(1)}\ra S_\alpha^{(2)}\ra S_\alpha^{(3)}\ra
S_\alpha^{(4)}=:S_{\alpha+1},
$$
where $S_\alpha^{(1)}$ is a minimal seed (and so a domain
by Proposition~\ref{mindom}),
$S_\alpha^{(2)}=(S_\alpha^{(1)})^+$ (an integral extension
of a seed and so a seed by Theorem~\ref{intseed}),
$S_\alpha^{(3)}$ is a quasilocal balanced big Cohen-Macaulay
$R$-algebra (which $S_\alpha^{(2)}$ maps to by
Lemma~\ref{local}), and $S_\alpha^{(4)}$ is the
$m$-adic completion of $S_\alpha^{(3)}$ (which is
$m$-adically separated and a balanced big Cohen-Macaulay
algebra by \cite[Theorem 1.7]{Bar-Str}).
If $\alpha$ is a limit ordinal, then we will define
$$
S_\alpha:=\varinjlim_{\gamma<\alpha} S_\gamma.
$$
Given our definition for $S_\alpha$ for each ordinal
$\alpha<\beta$, we define
$$
B:=S_{\beta}=\varinjlim_{\alpha<\beta} S_\alpha.
$$

Since each $S_\alpha$ maps to a domain $S_\alpha^{(1)}$,
and conversely, each $S_\alpha^{(1)}$ maps to
$S_{\alpha+1}$, the ring $B$ can be written as a direct
limit of domains and is, therefore, also a domain. Similarly,
$B$ is also a direct limit of absolutely integrally closed
domains (using $S_\alpha^{(2)}$ for each $\alpha<\beta$).
If we let $L$ be the algebraic closure of the fraction field
of $B$ and let $K_\alpha$ be the algebraic closure of the fraction
field of $S_\alpha^{(2)}$ for each $\alpha<\beta$, then an element
$u\in L$ satisfying a monic polynomial equation over $B$ is the image
of an element $v$ in some $K_\alpha$ that satisfies a monic
polynomial equation over $S_\alpha^{(2)}$ for some $\alpha$. Since
$S_\alpha^{(2)}$ is absolutely integrally closed, $v$ is in
$S_\alpha^{(2)}$, and so its image $u$ in $L$ is also in $B$. Therefore,
$B$ is absolutely integrally closed.

Using the rings
$S_\alpha^{(3)}$, we can see that $B$ is the direct limit
of quasilocal balanced big Cohen-Macaulay $R$-algebras, and so $B$
is itself a quasilocal balanced big Cohen-Macaulay algebra. Indeed, it is
easy to see that 
$$
(x_1,\ldots,x_k)B:_B x_{k+1}\subseteq (x_1,
\ldots x_k)B
$$ 
for each partial system of parameters
$x_1,\ldots,x_{k+1}$ of $R$ as this fact is true in each
$S_\alpha^{(3)}$. Furthermore, $mB\neq B$, as the
opposite would imply that $mS_\alpha^{(3)} = S_\alpha^{(3)}$
for some $\alpha$, which is impossible in a balanced big
Cohen-Macaulay algebra. It is also straightforward
to verify that a direct limit of quasilocal rings
is quasilocal.

Finally, to see that $B$ is $m$-adically separated, we note that
$B$ is a direct limit of the $m$-adically separated rings
$S_\alpha^{(4)}$. Suppose that
$u\in \bigcap_k m^kB$. Then for each $k$, there exists
an ordinal $\alpha(k)$ such that
$u\in m^k S_{\alpha(k)}^{(4)}$. Let $\alpha$ be the union
of all the $\alpha(k)$. Since $\alpha(k)<\beta$ for all
$k$ and $\beta$ is uncountable, we see that $\alpha<\beta$. Therefore,
$u\in \bigcap_k m^k S_\alpha$, and so
$u\in \bigcap_k m^k S_\alpha^{(4)} = 0$.
\end{proof}

\section{Tensor Products and Base Change}

In this section, we provide positive answers to two previously open questions about
big Cohen-Macaulay algebras.
Given two big Cohen-Macaulay $R$-algebras
$B$ and $B'$ over a complete local domain,
does there exist a big Cohen-Macaulay algebra $C$ such that
$$\xymatrix{
B\ar[r] & C \\
R\ar[r]\ar[u] & B'\ar[u] }$$
commutes?

Another open question involves base change
$R\ra S$ between complete local
domains. Given a big Cohen-Macaulay $R$-algebra $B$, can $B$ be mapped to a big
Cohen-Macaulay $S$-algebra $C$ such that the diagram
$$\xymatrix{
B\ar[r] & C \\
R\ar[r]\ar[u] & S\ar[u] }$$
commutes?

We will show that both of these questions have positive answers in
positive characteristic. First, we address the question of why the
tensor product of seeds is a seed. We will derive our result from
the case of a regular local base ring and make use of tight
closure and test elements for the general case.

By applying $T\tensor_S -$ to the sequence $0\ra I:_{S}x/I\ra
S/I\overset{x}\ra S/I$ one has:

\begin{lemma} If $S$ and $T$ are any commutative rings
such that $T$ is flat over $S$, $I$ is an ideal of $S$, and $x\in S$, then
$IT:_{T}x = (I:_{S}x)T$. \end{lemma}

\begin{lemma} If $C$ is a (balanced) big Cohen-Macaulay algebra over a local
ring $(S,n)$ and $D$ is faithfully flat over $C$, then $D$ is a (balanced) big
Cohen-Macaulay $S$-algebra. \end{lemma}
\begin{proof} Let $x_{1},\ldots,x_{k+1}$ be part of a system of parameters for
$S$, and let $d$ be an element of $(x_{1},\ldots,x_{k})D:_{D}x_{k+1}$.  Using
the lemma above,
$$d\in ((x_{1},\ldots,x_{k})C:_{C}x_{k+1})D\subseteq
((x_{1},\ldots,x_{k})C)D = (x_{1},\ldots,x_{k})D$$
because $C$ is a (balanced) big Cohen-Macaulay $S$-algebra. 
Finally, as $D$ is faithfully
flat over $C$, and $nC\neq C$, $nD\neq D$ either.
\end{proof}

If $A$ is a regular local ring, and $B$ and $B'$ are balanced 
big Cohen-Macaulay
$A$-algebras, then $B$ is faithfully flat over $A$
by Proposition~\ref{regflatbigCM}. Therefore, $B\tensor_{A}B'$
is faithfully flat over $B'$, and since $B'$ is a balanced big Cohen-Macaulay
$A$-algebra, we can use the previous lemma to conclude:

\begin{lemma}\label{prodreg} If $A$ is a regular local Noetherian
ring, and $B$ and $B'$ are balanced big Cohen-Macaulay $A$-algebras, 
then $B\tensor_{A}B'$
is a balanced big Cohen-Macaulay $A$-algebra as well. 
Consequently, if $S$ and $S'$ are
seeds over $A$, then $S\tensor_{A}S'$ is a seed over $A$.
\end{lemma}

We can now establish our first result, concerning tensor products of
seeds.

\begin{thm}\label{product} Let $(R,m)$ be a complete local domain
of positive characteristic. If $(S_i)_{i\in I}$ is an arbitrary family
of seeds over $R$, then
$\bigotimes_{i\in I} S_i$ is also a seed over $R$. Consequently, if $B$ and
$B'$ are (balanced) big Cohen-Macaulay $R$-algebras, 
then there exists a balanced big Cohen-Macaulay
$R$-algebra $C$ filling the commutative diagram:
$$\xymatrix{
B\ar[r] & C \\
R\ar[r]\ar[u] & B'\ar[u] }$$
\end{thm}
\begin{proof} Since a direct limit of seeds is a seed
by Lemma~\ref{limit}, we may assume that $I$ is a finite set. By induction,
we may assume that $I$ consists of two elements.
We may then also assume that $S_1=B$ and $S_2=B'$ are balanced big
Cohen-Macaulay $R$-algebras.

By the Cohen Structure Theorem, $R$ is a
module-finite extension of a complete regular local ring $A$.
We reduce to the case that $R$ is a separable extension of $A$.
For any $q=p^{e}$, $R[A^{1/q}]$ is still a module-finite
extension of $A^{1/q}$, a complete regular local ring. Furthermore,
$B[A^{1/q}]$ and $B'[A^{1/q}]$ are clearly still balanced
big Cohen-Macaulay $R[A^{1/q}]$-algebras, and
if $B[A^{1/q}]\tensor_{R[A^{1/q}]}B'[A^{1/q}]$ is a seed over
$R[A^{1/q}]$, then it is also a seed over $R$, which will show that
$B\tensor_{R}B'$ is a seed over $R$. We may therefore
replace $A$ and $R$ by $A^{1/q}$ and $R[A^{1/q}]$ for any $q\geq 1$.
We claim that for any $q\gg 1$, we obtain a separable extension
$A^{1/q}\ra R[A^{1/q}]$. Indeed, suppose that $R$ is not separable
over $A$. If $K$ is the fraction field of $A$, then
$K\tensor_{A}R$ is a finite product of finite field extensions of $K$,
one of which is not separable over $K$.
It suffices to let $L$ be a finite inseparable field
extension of $K$ and show that $L[K^{1/q}]$ is separable over $K^{1/q}$,
for some $q\gg 1$. For
any element $y$ of $L$ whose minimal polynomial is inseparable, we can
find $q$ sufficiently large so that the minimal polynomial of $y$ in
$L[K^{1/q}]$ over $K^{1/q}$ becomes separable. Since $L$ is
a finite extension, for any $q$ sufficiently large
$L[K^{1/q}]$ becomes separable over $K^{1/q}$, which implies
that $R[A^{1/q}]$ will be separable over $A^{1/q}$.

We can now assume without loss of generality that $R$ is separable
over $A$. Let $J$ be the ideal of $R\tensor_{A}R$ generated by all
elements killed by an element of $A$.
Then $R_{0}:= (R\tensor_{A}R)/J$
is a module-finite extension of $A$  and a
reduced ring as $R_{0}$ is also a separable extension of $A$.

Since $B$ and $B'$ are balanced big Cohen-Macaulay $R$-algebras and $R$ is a
module-finite extension of $A$, we also have that $B$ and $B'$ are balanced big
Cohen-Macaulay \mbox{$A$-algebras}. By Lemma~\ref{prodreg}, $B\tensor_{A}B'$ is also a
balanced big Cohen-Macaulay \mbox{$A$-algebra}. 
The $R$-algebra structures of $B$ and $B'$ induce
a natural map from $R\tensor_{A}R$ to $B\tensor_{A}B'$. Since the
latter ring is balanced big Cohen-Macaulay over $A$, the ideal $J$
is contained in the kernel of the map to $B\tensor_{A} B'$.
Therefore, this map factors through $R_{0}$.

Since $R$ is a domain and a homomorphic image of $R\tensor_{A}R$ (and
so of $R_{0}$), the kernel of $R_{0}\ra R$ is a
prime ideal $P$, and the
kernel of $B\tensor_{A}B'\ra B\tensor_{R}B'$ is the extended ideal
$P(B\tensor_{A}B')$. Since $R_{0}$ is a module-finite
extension of $A$, $\dim R=\dim R_{0}$, so
that $P$ is a minimal prime of $R_{0}$.
We therefore obtain a commutative diagram
$$\xymatrix{B\tensor_{A}B'\ar[r] & B\tensor_{R}B' \\
R_{0}\ar[r]\ar[u] & R\ar[u]}
$$
where the horizontal maps are the result of killing the minimal prime
ideal $P$ of $R_{0}$ (resp., $P(B\tensor_A B')$).

We next construct a durable colon-killer in $B\tensor_{R}B'$.
Since $R_{0}$ is reduced and $P$ is a
minimal prime, there exists $c'\not\in P$ such that $c'$ kills $P$
and $P(B\tensor_{A}B')$.
As $c'\not\in P$, its image in $R$ is nonzero. Since $R$ is a
complete local domain, $R$ has a test
element $c''\neq 0$ by Theorem~\ref{testexist}. Let $c=c'c''$, a
nonzero test element of $R$, and let $d$ be a lifting of $c$ to
$R_{0}$ so that $dP(B\tensor_{A}B') = 0$.

Now, suppose that $x_{1},\ldots,x_{n}$ is a system of parameters of
$A$, and suppose that
$x_{k+1}u\in(x_{1},\ldots,x_{k})(B\tensor_{R}B')$ for some $k\leq n-1$.
Then
$$x_{k+1}u\in ((x_{1},\ldots,x_{k}) + P)(B\tensor_{A}B'),$$
and $x_{k+1}du\in (x_{1},\ldots,x_{k})(B\tensor_{A}B')$. Since
$B\tensor_{A}B'$ is a balanced big Cohen-Macaulay $A$-algebra,
$du\in (x_{1},\ldots,x_{k})(B\tensor_{A}B')$,
and so $cu\in
(x_{1},\ldots,x_{k})(B\tensor_{R}B')$, as $c$ and $d$ have the same
image in $B\tensor_{R}B'$. Therefore, $c$ is a colon-killer for systems
of parameters of $A$ in $B\tensor_{R}B'$. By Proposition~\ref{ckbasechange},
we may replace $c$ by a power that is a colon-killer
for systems of parameters of $R$.

Finally, suppose that $c^{N}\in\bigcap_{k} m^{k}(B\tensor_{R}B')$. Then
\cite[Theorem 8.6]{Ho94} and Theorem~\ref{bigCMsolid} imply that
$c^{N}\in\bigcap_{k} ((m^{k})^{*})$ since $B\tensor_{R}B'$ is solid over
$R$ by Proposition~\ref{solidprop}(b). 
As $c$ is a test element in $R$, $c^{N+1}\in
\bigcap_{k}m^{k} = 0$, a contradiction. Hence, $c$ is a durable
colon-killer in $B\tensor_{R}B'$, and so $B\tensor_{R}B'$ is a seed
over $R$ by Theorem~\ref{ckseed}.
\end{proof}

We now proceed to the question of whether being a seed over a complete
local domain of positive characteristic is a property that is
preserved by base change to another complete local domain. If $R\ra
S$ is a map of complete local rings, then \cite[Theorem 1.1]{AFH} says
that the map factors through a complete local ring $R'$ such that
$R\ra R'$ is faithfully flat with regular closed fiber, and $R'\ra S$
is surjective. It therefore suffices to treat the cases of a flat local
map with regular closed fiber and the case of a surjective map
with kernel a prime ideal.
We start with an elementary lemma and
then prove the result for the flat local case.

\begin{lemma}\label{flatseq}
Let $(A,m)$ be a local Noetherian ring, and let $B$ be
a flat, local Noetherian $R$-algebra. If $y_{1},\ldots,y_{t}$ is a
regular sequence on $B/mB$, then $y_{1},\ldots,y_{t}$ is a regular
sequence on $B$, and $B/(y_{1},\ldots,y_{t})B$ is flat over $A$.
\end{lemma}
\begin{proof} The proof is immediate by induction on $t$, where the base case
of $t=1$ is given by \cite[(20.F)]{Mat}.
\end{proof}

With this lemma, we are ready to prove our base change result
for flat local maps. It is perhaps
interesting to note that our argument only requires that
the closed fiber is Cohen-Macaulay, not regular. Unlike the surjective case,
we will not need to assume that our rings have positive
characteristic, are complete, or are domains.

\begin{prop}\label{seedflat} Let $R\ra S$ be a flat local map
of Noetherian local rings with a Cohen-Macaulay closed fiber $S/mS$,
where $m$ is the maximal ideal of $R$. If $T$ is a seed over
$R$, then $T\tensor_R S$ is a seed over $S$.
\end{prop}
\begin{proof} It suffices to assume that $T=B$ is
a balanced big Cohen-Macaulay \mbox{$R$-algebra} and 
show that a single system of parameters of $S$ is a regular
sequence on $B\tensor_R S$.

Fix a system of parameters $x_1,\ldots,x_d$ for $R$. Since $S$ is
faithfully flat over $R$, we have the dimension equality
$\dim S = \dim R + \dim S/mS$. Hence, the images of
$x_1,\ldots,x_d$ in $S$ can be extended to a full system of parameters
$x_1,\ldots,x_d,\ldots, x_n$ of $S$, where
$x_{d+1},\ldots,x_{n}$ is a system of parameters for $S/mS$. As
$x_1,\ldots,x_d$ form a regular sequence on $B$  for
any $1\leq k\leq d-1$, we have an exact sequence:
$$
0\ra B/(x_1,\ldots,x_k)B \overset{x_{k+1}}\longrightarrow
B/(x_1,\ldots,x_k)B,
$$
and since $S$ is flat over $R$, the sequence remains
exact after tensoring with $S$. Thus,
$x_{k+1}$ is not a zerodivisor on
$$
(B/(x_1,\ldots,x_k)B)\tensor_R S\isom
(B\tensor_R S)/(x_1,\ldots,x_k)(B\tensor_R S)
$$
for any $1\leq k\leq d-1$.

It now suffices to show
that $x_{d+1},\ldots,x_n$ is a regular sequence on
the quotient
$$\overline{B}:=(B\tensor_R S)/(x_1,\ldots,x_k)(B\tensor_R S).$$
Let $I:=(x_{1},\ldots,x_{d})R$, and let $\overline{R}:=R/I$, and
$\overline{S}:=S/IS$. Then $\overline{S}$ is faithfully flat over
$\overline{R}$, and since $x_{d+1},\ldots,x_{n}$ is a system of parameters for
the Cohen-Macaulay ring $S/mS$, it is a regular sequence on $S/mS\isom
\overline{S}/\overline{m}\overline{S}$, where $\overline{m}=m/I$, the
maximal ideal of $\overline{R}$.
We can now apply Lemma~\ref{flatseq} to $\overline{R}$ and
$\overline{S}$ to conclude that $x_{d+1},\ldots,x_{n}$ is a regular
sequence on $\overline{S}$ and that
$\overline{S}/(x_{d+1},\ldots,x_{k})\overline{S}$ is flat over
$\overline{R}$ for all $d+1\leq k\leq n$.

For any $d\leq k\leq n-1$, we have a short exact sequence
$$
0\ra \overline{S}/(x_{d+1},\ldots,x_{k})\overline{S}
\overset{x_{k+1}}\longrightarrow
\overline{S}/(x_{d+1},\ldots,x_{k})\overline{S} \ra
\overline{S}/(x_{d+1},\ldots,x_{k},x_{k+1})\overline{S} \ra 0,
$$
where $x_{d+1},\ldots,x_{k}$ is the empty sequence when $k=d$.
Since $\overline{S}/(x_{d+1},\ldots,x_{k+1})\overline{S}$ is flat
over $\overline{R}$, we have
$\Tor_{1}^{\overline{R}}(\overline{B},
\overline{S}/(x_{d+1},\ldots,x_{k+1})\overline{S}) = 0.$
Therefore,
$$
0\ra \overline{B}\tensor_{\overline{R}}
(\overline{S}/(x_{d+1},\ldots,x_{k})\overline{S})
\overset{x_{k+1}}\longrightarrow
\overline{B}\tensor_{\overline{R}}
(\overline{S}/(x_{d+1},\ldots,x_{k})\overline{S})
$$
is exact, and since
$$
\overline{B}\tensor_{\overline{R}}
(\overline{S}/(x_{d+1},\ldots,x_{k})\overline{S}) \isom
(\overline{B}\tensor_{\overline{R}} \overline{S})/
(x_{d+1},\ldots,x_{k})
(\overline{B}\tensor_{\overline{R}} \overline{S}),
$$
$x_{d+1},\ldots,x_{n}$ is a possibly improper regular sequence on
$\overline{B}\tensor_{\overline{R}} \overline{S}$.
We can now finally see that $x_{1},\ldots,x_{d},x_{d+1},\ldots,x_{n}$ is a
possibly improper regular sequence on $B\tensor_{R}S$ as
$
\overline{B}\tensor_{\overline{R}} \overline{S}\isom
(B\tensor_{R}S)/I(B\tensor_{R}S).
$
If, however, $(B\tensor_{R}S)/n(B\tensor_{R}S) = 0$, where
$n$ is the maximal ideal of $S$, then $B\tensor_{R}(S/nS)=0$, which
implies that the product $(B/mB)\tensor_{R/m}(S/nS)=0$ over the
field $R/m$. Therefore, $B/mB=0$ or $S/nS=0$, but neither of these occurs, so
we have a contradiction. Hence, $x_{1},\ldots,x_{n}$ is a regular
sequence on $B\tensor_{R}S$, and so $B\tensor_{R}S$ is a seed.
\end{proof}

We will now handle the case of a surjective map $R\ra S$ of
complete local domains of positive characteristic.
First we demonstrate the result when $R$ is normal and then
show how the problem can be reduced to the normal case using 
Theorem~\ref{intseed}.

In the normal
case, we will make use of test elements again. When $R$ is normal,
the defining ideal $I$ of the singular locus has height at least 2. Since $R_{c}$ is
regular for any $c\in I$,
\cite[Theorem 6.1]{HH94sm} implies that some power $c^{N}$ is a test
element when $R$ is a reduced excellent local
ring. Hence, if $P$ is a height 1 prime of $R$, there exists
a test element $c$ of $R$ not in $P$. We record this fact in the
following lemma.

\begin{lemma}\label{normaltest}
Let $R$ be a normal excellent local ring of positive characteristic.
If $P$ is a height 1 prime of $R$, then there exists a test
element $c\in R\setminus P$.
\end{lemma}

\begin{lemma}\label{normalcase}
Let $(R,m)$ be a complete local, normal domain of positive
characteristic, and let $S=R/P$,
where $P$ is a height 1 prime of $R$. If $T$ is a seed over $R$,
then $T/PT$ is a seed over $S$.
\end{lemma}
\begin{proof} We can assume that $B=T$ is a 
balanced big Cohen-Macaulay \mbox{$R$-algebra}.
Since $R$ is normal and $\hgt\ P=1$, we see that $R_P$ is a DVR.
Therefore, $PR_P$ is a principal ideal generated
by the image of $x\in R$. Each element of $P$ is then multiplied
into $xR$ by some element of $R\setminus P$, and since $P$ is finitely
generated, there exists $c'\in R\setminus P$ such that $c'P\subseteq
xR$. By Lemma~\ref{normaltest}, there exists a test element
$c''\in R\setminus P$, and so if we put $c:=c'c''$, then $c$ is a test
element, $cP\subseteq xR$, and $c$ is not in $P$.

We claim that $c$ is a weak durable colon-killer for $B/PB$.
Extend $x$ to a full system of parameters $x,x_2,\ldots,x_d$ for $R$. Then
$x_2^t,\ldots,x_d^t$ is a system of parameters for $S$ for any $t\in\nn$.
Suppose that $bx_{k+1}^t\in(x_2^t,\ldots,x_k^t)B/PB$ for some
$k\leq d-1$ and some $t$. This relation lifts to a relation
$bx_{k+1}^t\in(x_2^t,\ldots,x_k^t)B + PB$ in $B$, and so
$cbx_{k+1}^t\in (x,x_2^t,\ldots,x_k^t)B$.
Since $B$ is a balanced big Cohen-Macaulay $R$-algebra, we have
$cb\in (x,x_2^t,\ldots,x_k^t)B$, and so $cb\in (x_2^t,\ldots,x_k^t)
B/PB$.

To finish, suppose that $c^N\in\bigcap_k (m/P)^kB/PB$.
We can then lift to $B$ to obtain
$c^N\in \bigcap_k (m^k+P)B$. Since $R$ is a complete local domain,
and $B$ is a balanced big Cohen-Macaulay $R$-algebra, 
Theorem~\ref{bigCMTC}
implies that $c^N\in\bigcap_k (m^k+P)^*$. Since $c$ was
chosen to be a test element, we have
$c^{N+1}\in\bigcap_k (m^k+P) = P$, a contradiction as 
$c\not\in P$. Therefore,
the image of $c$ in $B/PB$ is a weak durable colon-killer,
and so $B/PB$ is a seed over $S$ by Theorem~\ref{ckseed}.
\end{proof}

We finally treat the case of an arbitrary surjection of
complete local domains by reducing to the case of
the previous lemma.

\begin{prop}\label{seedsurj} Let $R\ra S$ be
a surjective map of positive characteristic
complete local domains.
If $T$ is a seed over $R$, then $T\tensor_R S$
is a seed over $S$.
\end{prop}
\begin{proof} We can immediately assume that the kernel of
$R\ra S$ is a height 1 prime $P$ of $R$.
Let $R'$ be the normalization of $R$ in its fraction field.
Then $R'$ is also a complete local domain.
Since $R'$ is an integral extension of $R$, there exists
a height 1 prime $Q$ lying over $P$.

By Corollary~\ref{dom},
$T$ maps to a big Cohen-Macaulay $R$-algebra domain $B$, and so
we may replace $T$ by $B$ and assume that
$T$ is a domain. We then have
an integral extension $T[R']$ of $T$ inside the fraction field
of $T$. Since $T$ is a seed over $R$, Theorem~\ref{intseed}
implies that $T[R']$ is also a seed over $R$. Therefore,
$T[R']$ maps to a balanced big Cohen-Macaulay $R$-algebra $C$ (which is
also an $R'$-algebra), and so Corollary~\ref{bigCMbasechange}
implies that $C$
is a balanced big Cohen-Macaulay $R'$-algebra.
We now have the commutative diagram:
$$
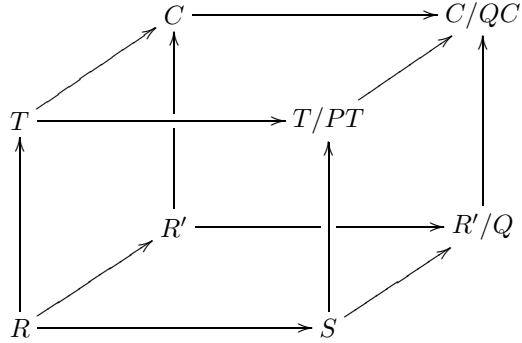
\xymatrix@!C{ & C\ar[rr] && C/QC \\
T\ar[ur]\ar[rr] & \ar[u] & T/PT\ar[ur] & \\
& R'\ar@{-}[r]\ar@{-}[u] & \ar[r] & R'/Q\ar[uu] \\
R\ar[ur]\ar[rr]\ar[uu] && S\ar[uu]\ar[ur] &  }
$$
since $S=R/P$ and $Q$ lies over $P$. By Lemma~\ref{normalcase}, 
$C/QC$ is a seed over $R'/Q$. Since
$R'/Q$ is an integral extension of $S=R/P$, every system of parameters of
$S$ is a system of parameters of $R'/Q$, and so $C/QC$ is
also a seed over $S$, but this implies that
$T/PT$ is also a seed over $S$, as needed.
\end{proof}

Now that we have shown that the property of being a seed
over a complete local domain is preserved by flat base change
with a regular closed fiber (Proposition~\ref{seedflat}) and by
surjections (Proposition~\ref{seedsurj}), we may apply
\cite[Theorem 1.1]{AFH} to factor any map of complete
local domains into these two maps. We therefore arrive
at the following theorem, which answers the base change
question asked at the beginning of the section.

\begin{thm}\label{basechange} Let $R\ra S$
be a local map of positive characteristic complete
local domains. If $T$ is a seed over $R$, then
$T\tensor_R S$ is a seed over $S$. Consequently, if $B$
is a big Cohen-Macaulay $R$-algebra, then there exists a 
balanced big Cohen-Macaulay
$S$-algebra $C$ filling the commutative diagram:
$$
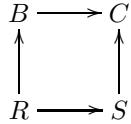
\xymatrix{
B\ar[r] & C \\
R\ar[r]\ar[u] & S\ar[u] }$$
\end{thm}

\section{Seeds and Tight Closure in Positive Characteristic}

Using Theorems~\ref{product} and \ref{basechange},
we can now use the class of all balanced big Cohen-Macaulay $R$-algebras
$\sB(R)$, where $R$ is a complete local domain of characteristic $p$,
to define a closure operation for all Noetherian rings of positive
characteristic. A key point is that
$\sB(R)$ is a directed family and has the right base change
properties. (It is easy to see that the two definitions of
$\natural$-closure below coincide for complete local domains.)

\begin{defn} Let $R$ be a complete local domain of positive characteristic,
and let $N\subseteq M$ be finitely
generated $R$-modules. Let $N_M^\natural$ be the set of all
elements $u\in M$ such that
$1\tensor u\in \im(B\tensor_{R} N\ra B\tensor_{R} M)$ for some
big Cohen-Macaulay $R$-algebra $B$.

Let $S$ be a Noetherian ring of positive characteristic,
and let $N\subseteq M$ be finitely generated $S$-modules.
Let $N_M^\natural$ be the set of all $u\in M$ such that for all
$S$-algebras $T$, where $T$ is a complete local domain, $1\tensor u
\in\im(T\tensor_{S} N\ra T\tensor_{S} M)_{T\tensor_{S} M}^\natural$.
We will call $N_{M}^{\natural}$ the \textit{$\natural$-closure}
of $N$ in $M$. \end{defn}

By Theorem~\ref{bigCMTC},
our new closure operation is equivalent to tight closure for complete
local domains of positive characteristic, but the
results above imply many of the properties one
would want in a good closure operation directly from the properties of
big Cohen-Macaulay algebras, independent of tight closure.
These properties include
persistence, $(IS)^{\natural}_{S}\cap R\subseteq I_{R}^{\natural}$ for
module-finite extensions, $I^{\natural}=I$ for ideals in a regular
ring, phantom acyclicity, and colon-capturing. (See Theorem~\ref{TCprop}
and \cite[Section 9]{HH90} for the statement of these properties for tight
closure. See \cite{D} for the statements and proofs for $\natural$-closure.)

These results add evidence to
the idea that such a closure operation can be defined for more general
classes of rings.

\section*{Acknowledgments}
I thank my advisor, Mel Hochster, for all of the time, help, and insight
given to me during the production of this work. I also thank the
entire Mathematics Department at the University of Michigan for
the five years of education and support provided to me while
I was a graduate student.


\begin{thebibliography}{HHH}

\bibitem[AFH]{AFH} \textsc{L. Avramov}, \textsc{H.-B. Foxby}, and \textsc{B.
Herzog}, \textit{Structure of local homomorphisms},
J. Algebra \textbf{164} (1994), no. 1, 124--145.

\bibitem[BS]{Bar-Str} \textsc{J. Bartijn} and \textsc{J.R. Strooker},
Modifications Monomiales, \textit{S\'{e}minaire d'Alg\`{e}bre
Dubreil-Malliavin, Paris 1982},
Lecture Notes in Mathematics, \textbf{1029},
Springer-Verlag, Berlin, 1983, 192--217.

\bibitem[BH]{BH} \textsc{W. Bruns} and \textsc{J. Herzog}, \textit{Cohen-Macaulay
Rings}, Cambridge Studies in Advanced Mathematics \textbf{39},
Cambridge University Press, Cambridge, 1993.

\bibitem[D]{D} \textsc{G. Dietz}, \textit{Closure Operations in Positive Characteristic
and Big Cohen-Macaulay Algebras}, Thesis, Univ. of Michigan, 2005.

\bibitem[E]{E} \textsc{D. Eisenbud}, \textit{Commutative Algebra with a View
Toward Algebraic Geometry}, Graduate Texts in Mathematics \textbf{150},
Springer-Verlag, New York, 1995.

\bibitem[Heit]{Heit02} \textsc{R. Heitmann}, \textit{The direct summand conjecture 
in dimension three}, Ann. of Math. (2) \textbf{156} (2002), no. 2, 695--712.

\bibitem[Ho1]{Ho75}  \textsc{M. Hochster}, \textit{Topics in the homological theory
of modules over commutative rings}, CBMS Regional Conference Series \textbf{24},
Amer. Math. Soc., 1975.

\bibitem[Ho2]{Ho94} \textsc{M. Hochster}, \textit{Solid closure}, Contemp. Math.
\textbf{159} (1994), 103--172.

\bibitem[Ho3]{Ho02} \textsc{M. Hochster}, \textit{Big Cohen-Macaulay algebras in 
dimension three via Heitmann's theorem}, J. Algebra \textbf{254} (2002), no. 2, 395--408.

\bibitem[HH1]{HH90} \textsc{M. Hochster} and \textsc{C. Huneke}, \textit{Tight
closure, invariant theory, and the Brian\c{c}on-Skoda theorem}, J. Amer.
Math. Soc. \textbf{3} (1990), 31--116.

\bibitem[HH2]{HH92} \textsc{M. Hochster} and \textsc{C. Huneke}, \textit{Infinite
integral extensions and big Cohen-Macaulay algebras}, Annals of Math.
\textbf{135} (1992), 53--89.

\bibitem[HH3]{HH94} \textsc{M. Hochster} and \textsc{C. Huneke}, \textit{Tight
closure of parameter ideals and splitting in module-finite  extensions}, J.
Algebraic Geom.  \textbf{3}  (1994),  no. 4, 599--670.

\bibitem[HH4]{HH94sm} \textsc{M. Hochster} and \textsc{C. Huneke},
\textit{$F$-regularity, test elements, and smooth base change},
Trans. Amer. Math. Soc. \textbf{346} (1994), no. 1, 1--62.

\bibitem[HH5]{HH95} \textsc{M. Hochster} and \textsc{C. Huneke},
\textit{Applications of the existence of big Cohen-Macaulay algebras},
Adv. Math.  \textbf{113}  (1995),  no. 1, 45--117.

\bibitem[HJ]{HJ} \textsc{K. Hrbacek} and \textsc{T. Jech},
\textit{Introduction
to Set Theory}, Second edition.
Monographs and Textbooks in Pure and Applied Mathematics, \textbf{85},
Marcel Dekker, Inc., New York, 1984.

\bibitem[Mat]{Mat} \textsc{H. Matsumura}, \textit{Commutative Algebra},
Benjamin-Cummings, New York, 1980.

\bibitem[N1]{N65} \textsc{M. Nagata}, \textit{Lectures on the Fourteenth Problem
of Hilbert}, Tata Institute of Fundamental Research, Bombay, 1965.

\bibitem[N2]{N72} \textsc{M. Nagata}, \textit{Local Rings}, Interscience,
New York, 1972.

\end{thebibliography}
\end{document}